\documentclass[11pt,a4paper]{amsart}

\usepackage[colorlinks=true,linkcolor=red!70!black,urlcolor=MidnightBlue,citecolor=MidnightBlue]{hyperref}

\usepackage[backend=bibtex, backref=true, giveninits=true, style=alphabetic,doi=false,bibencoding=utf8,language=auto,autolang=other,maxnames=10, maxalphanames=10]{biblatex}
\usepackage[margin=3cm,top=2.5cm,bottom=2cm]{geometry}

\usepackage{amsmath, amssymb, amsthm, mathrsfs}
\usepackage{graphicx,enumitem,microtype} 
\usepackage[dvipsnames]{xcolor}
\usepackage{thmtools}
\usepackage[nameinlink,capitalise]{cleveref}
\crefname{equation}{}{}

\usepackage[utf8]{inputenc} 
\usepackage[T1]{fontenc}
\usepackage{lmodern}

\newcommand{\ds}{\displaystyle}

\newcommand{\N}{\mathbb N}
\newcommand\R{{\mathbb R}}

\newcommand{\I}{\mathbb I}

\renewcommand{\S}{\mathbb{S}}

\newcommand\CC{\mathcal C}

\newcommand\FF{\mathcal F}

\newcommand\HH{\mathcal H}

\newcommand\KK{\mathcal K}
\newcommand\LL{\mathcal L}

\newcommand\PP{\mathcal P}

\newcommand\TT{\mathcal T}

\newcommand\WW{\mathcal W}

\newcommand\XX{\mathcal X}
\newcommand\YY{\mathcal Y}
\newcommand\ZZ{\mathcal Z}

\newcommand\eps{\varepsilon}

\newcommand{\indicatrice}{\chi}

\newcommand{\mean}{\mathsf{m}}

\newcommand{\ent}{\mathsf{Ent}}

\newcommand{\bra}[1]{\left( #1 \right)}
\newcommand{\sqa}[1]{\left[ #1 \right]}
\newcommand{\cur}[1]{\left\{ #1 \right\}}
\newcommand{\ang}[1]{\left< #1 \right>}
\newcommand{\abs}[1]{\left| #1 \right|}
\newcommand{\nor}[1]{\left\| #1 \right\|}

\newcommand{\Var}{\operatorname{Var}}
\newcommand{\osc}{\operatorname{osc}}
\renewcommand{\I}{\mathcal{Z}}

\newcommand{\la}{\langle}
\newcommand{\ra}{\rangle}

\newcommand{\Lip}{\operatorname{Lip}}
\newcommand{\diam}{\operatorname{diam}}

\numberwithin{equation}{section}

\setlist[enumerate]{wide,labelindent=0cm,label=\textnormal{(\arabic*)},itemsep=5pt,topsep=4pt}

\theoremstyle{plain}
\newtheorem{theo}{Theorem}[section]
\newtheorem{prop}[theo]{Proposition}
\newtheorem{lem}[theo]{Lemma}

\theoremstyle{remark}
\newtheorem{rem}[theo]{Remark}
\newtheorem{ex}[theo]{Example}

\theoremstyle{definition}
\newtheorem{defi}[theo]{Definition}

\title[Stability in optimal transport for power costs]{Quantitative stability in optimal transport for general power costs}

\author[O. Mischler]{Octave Mischler}
\address[Octave Mischler]{DMA,
  Ecole Normale Supérieure PSL, 45 rue d'Ulm,
  75005 Paris 5, France}
\email{octave.mischler@ens.psl.eu}
 
 \author[D. Trevisan]{Dario Trevisan}
\address[Dario Trevisan]{Dipartimento di Matematica, Università di Pisa,
Largo Bruno Pontecorvo 5,
56127 - Pisa (PI) Italy}
\email{dario.trevisan@unipi.it}

\date{\today}


\begin{document}

\begin{abstract}
We establish novel quantitative stability results for optimal transport problems with respect to perturbations in the target measure. We provide bounds on the stability of optimal transport potentials and maps, which are relevant for both theoretical and practical applications. Compared to previous results, ours apply to a wide range of costs, including all Wasserstein distances with power cost exponent strictly larger than $1$ and leverage mostly assumptions on the source measure, such as log-concavity and bounded support. Our proofs follow the same strategy as [Delalande, A., {\it Quant. Stab. in Quad. Opt. Trans., PhD thesis, 2022}] up to several technical improvements. Our work provides a significant step forward in the understanding of stability of optimal transport problems, as previous results were mostly limited to the case of the quadratic cost.
\end{abstract}

\maketitle


\section{Introduction} 

 A well-posed mathematical problem in the Hadamard sense is one for which existence and uniqueness of solutions hold, but also stability with respect to perturbations in the data. Although qualitative stability in most situations emerges as a by-product of existence and uniqueness arguments, quantitative stability becomes particularly crucial in applications, for it directly relates to the convergence of numerical methods. This applies in particular to optimal transport problems, which in the last decades have had significant impact across various fields of applied mathematics, most notably PDEs of evolution type \cite{AGS} and machine learning \cite{peyre2019computational}. It is also relevant in statistics, when one considers plug-in estimators from i.i.d. samples. In this work, we establish novel quantitative stability results for optimal transport problems with respect to a large variety of costs, including Wasserstein distances. 

 The Wasserstein distance of order $p \ge 1$ between two probability measures $\rho$, $\mu$ on $\R^d$, with finite $p$-th moment, is defined as the quantity
\begin{equation}
 \WW_p(\rho, \mu) := \inf_{\pi \in \Pi(\rho, \mu)} \left( \int_{\R^d \times \R^d} |y-x|^p d \pi(x, y) \right)^{1/p},
\end{equation}
where $\Pi(\rho, \mu)$ denotes the set of couplings between $\rho$ and $\mu$, i.e., (joint) probability measures on $\R^d\times \R^d$ with marginal distributions $\rho$, $\mu$, and $|\cdot|$ denotes the Euclidean norm. Under mild assumptions, in particular if $\rho$ is absolutely continuous with respect to the Lebesgue measure and $p>1$, one can argue \cite[Theorem 6.2.4]{AGS} that the optimizer is unique and given by an optimal transport map $T$ for $\WW_p(\rho, \mu)$, i.e., one has $y = T(x)$ for $\pi$-a.e.\  $(x,y)$ (with necessarily $T_\sharp \rho = \mu$) and therefore
\begin{equation}\label{eq:wp-map}
 \WW_p(\rho, \mu)^p =   \int_{\R^d} |T(x)-x|^p d \rho(x).
\end{equation}
For the quadratic cost case $p=2$, the optimal transport map is often referred to as Brenier's map.

Even when an optimal transport map exists, it is useful to define $\WW_p$ as a minimum over couplings, as one can exploit convex duality and obtain the equivalent formulation
\begin{equation}\label{eq:wp-dual}
  \WW_p(\rho, \mu)^p = \sup_{(\phi,\psi)} \cur{  \int_{\R^d} \phi d \rho + \int_{\R^d} \psi d \mu  },
\end{equation}
where the supremum runs among upper semi-continuous  functions $(\phi,\psi)$ such that $\phi(x)+\psi(y) \le |x-y|^p$ for every $x$, $y\in \R^d$.
An optimizing pair $(\phi,\psi)$ is called an optimal (or Kantorovich) pair of potentials for $\WW_p(\rho, \mu)$. Under minimal assumptions, one can always argue that an optimizer exists, together with a qualitative stability result stating convergence of the potentials or the maps, see e.g. \cite[Theorem 5.20, Corollary 5.23]{villani2009optimal}. 

\subsection{Main results}

The aim of this work is to establish new quantitative stability results for optimal transport potentials and maps with respect to perturbations of the ``target'' measure $\mu$, keeping one ``source'' measure  $\rho$ fixed. In a slightly simplified form, to ease the notation, our main contributions can be summarized in the following two results. They extend and improve results from \cite{delalande2022quantitative} (quadratic case only) to a wider range of costs and more general assumptions on the measure by following the same line of proof but overcoming several technical difficulties along the way. The first one concerns stability of Kantorovich potentials and applies to any $\rho$ with bounded density with respect to a $\log$-concave probability measure with bounded support.

\begin{theo}\label{thm:main-pot}
Let $d \ge1$, $\lambda$ be a $\log$-concave probability measure on $\R^d$ with bounded support,  $\rho$ be a probability measure on $\R^d$ absolutely continuous with respect to $\lambda$, with density bounded from above and below by strictly positive constants, let $\YY \subseteq \R^d$ be compact and $p>1$. Then, there exists $C = C(\lambda, \rho, \YY, p)< \infty$ such that, for any $\mu$, $\nu$, probability measures supported on $\YY$, it holds
\begin{equation}\label{eq:stability-potential}
 \| \phi_{\mu} - \phi_\nu \|_{L^2(\rho)} \le C \WW_1(\mu, \nu)^{\theta},
\end{equation}
with
\begin{equation}
\theta  =
\begin{cases}  1- \frac{1}{p} & \text{if $1<p<2$,}\\
\frac{ 1}{2} & \text{if $p \ge 2$,}
\end{cases}
\end{equation}
and where $(\phi_\mu, \psi_\mu)$ denotes the Kantorovich potentials for $\WW_p(\rho, \mu)$, i.e., the optimizer in \cref{eq:wp-dual}, with zero $\rho$ mean, i.e., $\int \phi_{\mu} d \rho = 0$ (and similarly $(\phi_\nu, \psi_\nu)$, for $\WW_p(\rho, \nu)$).
\end{theo}

 Note that it is classically known that the Kantorovich potential with zero $\rho$-mean is unique in this setting, which can also be recovered by taking $\mu = \nu$. The only assumption on $\lambda$ being $\log$-concavity, it allows measures with densities going to $0$ at the boundary of their support. Our second main result deals instead with optimal transport maps: in this case, we require $\rho$ to have a density bounded from above and below with respect to the Lebesgue measure on a bounded convex set.

\begin{theo}\label{thm:main-map}
Let $d\ge 1$, $\rho$ be a probability measure on $\R^d$ with bounded and convex support, absolutely continuous with respect to the Lebesgue measure, with density bounded from above and below by strictly positive constants, let $\YY \subseteq \R^d$ be compact and $p>1$. Then,
\begin{equation}
\begin{cases} \text{for every $\theta \in \bra{ 0, \frac{ (p-1)^2}{p(p+1)} }$} & \text{if $1<p<2$,}\\
 \text{for $\theta = \frac{ 1}{6(p-1)}$} & \text{if $p \ge 2$,}
\end{cases}
\end{equation}
there exists $C = C(\rho, \YY, p,\theta)< \infty$ such that, for any $\mu$, $\nu$, probability measures supported on $\YY$, it holds
\begin{equation}\label{eq:stability-map}
 \| T_{\mu} - T_\nu \|_{L^2(\rho)} \le C \WW_1(\mu, \nu)^{\theta},
\end{equation}
where $T_\mu$ denotes the optimal transport map for $\WW_p(\rho, \mu)$, i.e.,  \cref{eq:wp-map} holds (and similarly $T_\nu$, for $\WW_p(\rho, \nu)$).
\end{theo}

\begin{rem}
It seems possible to write explicitly the constants $C=C(\rho, \YY, p, \theta)$ by carefully tracking their dependence in our arguments. Here, we use a simple scaling argument : in the case of the support of $\rho$ being a ball of radius $R$, with $\rho_{\min} R^{-d} \le \rho(x) \le \rho_{\max} R^{-d}$ a.e. in such ball, and $\mu$, $\nu$ contained also in the same ball, one finds
$$ \nor{ \phi_\mu - \phi_\nu}_{L^2(\rho)} \le C R^{p-\theta} \WW_1(\mu, \nu)^{\theta}$$
for some $C = C(\rho_{\min}, \rho_{\max}, p, d) < \infty$, and
$$ \nor{ T_\mu - T_\nu}_{L^2(\rho)} \le C R^{1-\theta} \WW_1(\mu, \nu)^{\theta}$$
for some $C = C(\rho_{\min}, \rho_{\max}, p, d, \theta) < \infty$. For more general (but still convex) supports, the dependence of the constants may be more involved, in particular in the case of the stability inequality for maps, because of our use of \Cref{prop:convgaliardo}, where the perimeter of the support explicitly appears. Note finally that in the inequality for potentials and when $p \geq 2$ the constant is explicit in \cref{theo:stabp>2}. \\
\end{rem}

\begin{rem}
Using the example of \cite[Lemma 2.2]{delalande2020} for any $p$, we can show that on a bounded set, the largest exponent $\theta$ for which \cref{eq:stability-map} holds is smaller than $\frac12$. However the ratio between this upper bound and our lower bound on the largest $\theta$ for which \cref{eq:stability-map} holds diverges as $p \to 1$ or $p \to \infty$, which leaves open the question of sharpness of our exponents. 
\end{rem}

\subsection{Related literature}

To our knowledge, most of the existing quantitative results concern quadratic optimal transport only. The first of these results is due to Ambrosio, reported by Gigli in \cite{Gigli2011} and reformulated in form closer to our setting in \cite{delalande2020}. It states stability of Brenier's map around measures $\mu$ such that $T_\mu$ is Lipschitz.

\begin{theo} [Ambrosio and Gigli]
Let $\rho$ be a probability density with bounded support $\XX \subseteq \R^d$, let $\YY \subset \R^d$ be  compact, and given $\mu, \nu \in \PP(\YY)$, write $T_{\mu}$ for the Brenier map from $\rho$ to $\mu$, i.e., \cref{eq:wp-map} holds with $p=2$ and $T=T_{\mu}$, similarly write $T_{\nu}$. Then,
\begin{equation}
\| T_{\mu} -T_{\nu}\|_{L^2(\rho)} \leq 2 \Lip(T_\mu) \diam(\XX) \WW_1(\mu,\nu)^{1/2},
\end{equation}
where $\Lip(\cdot)$ denotes the Lipschitz constant and $\diam(\cdot)$ the diameter.
\end{theo}

When compared with \cref{thm:main-map}, we see that the assumptions on $\rho$ are weaker, but are compensated by the fact that $T_{\mu}$ must be Lipschitz continuous, otherwise the bound is vacuous. In applications, regularity for the Brenier map is usually obtained as a consequence of Caffarelli's theory \cite{caffarelli1992regularity}, which requires stronger assumptions than \cref{thm:main-map}. As a recent extension to more general costs of this type of results, we mention \cite{gallouet2022strong}.

Other results in this direction have been obtained, notably \cite{neilan2018rates, berman2021convergence, li2021quantitative}, but to keep this part short, we only state one other result, which summarizes the results \cite[Theorem 5.14, Theorem 5.12, Corollary 5.8]{delalande2022quantitative}, first exposed in \cite{delalande2020, delalande2023quantitative}.

\begin{theo}[M\'erigot, Delalande and Chazal] \label{thm:delalande} Let $\XX \subseteq \R^d$ be a compact convex set, $\rho$ be a probability density supported on $\XX$, bounded from above and below by strictly positive constants. Let $p > \max\cur{d, 3}$ and assume that $\mu, \nu \in \PP(\R^d)$ have bounded $p$-th moment,
\begin{equation}
  \int_{\R^d} |x|^p d \mu(x) + \int_{\R^d }|x|^p d \nu(x) =: m_p < \infty.
\end{equation}
Then, there exists $C = C(d,p,\XX,\rho,m_p) < \infty$ such that
\begin{equation}\| \phi_\mu - \phi_\nu \|_{L^2(\rho)} \leq C \WW_1(\mu,\nu)^{1/2} \quad \text{and} \quad \|T_\mu -T_\nu \|_{L^2(\rho,\R^d)} \leq C \WW_1(\mu,\nu)^{1/(6+16d/p)},
\end{equation}
where $\phi_\mu$, $T_\mu$ denote respectively the Kantorovich potential with zero $\rho$-mean and the optimal transport map for $\WW_2(\rho, \mu)$ (and similarly $\phi_\nu$, $T_\nu$ for $\WW_2(\rho, \nu)$).

If moreover $\mu$,$\nu$ are supported on $\YY$ compact one can replace the exponent $1/(6+16d/p)$ above with $1/6$, for some constant $C = C(d,\XX,\rho, \YY)< \infty$.
\end{theo}

Let us finally mention \cite{divol2024tight} which provides stronger estimates than the previous result in the semi-discrete setting, and \cite{chizat2024sharper}, where improved rates for the Sinkhorn iterates are established using also arguments based on Prékopa-Leindler inequality, for general semi-concave ground cost, hence partly covering the $p\geq2$ part of the present paper. 

Technically, our work  improves upon the above results by extending the quantitative stability for optimal transport to larger families of costs, paving the way for a wider range of application. Our analysis, in particular \cref{thm:main-pot}, allows for a source probability $\rho$ that is merely $\log$-concave, and actually this assumption is relaxed by allowing for a density with respect to a given $\log$-concave measure  that is bounded from above and below by strictly positive constants; in a similar way,  in \cref{thm:main-map} the bounded support assumption for the target measures could be replaced by moment bounds -- see also extensions and open problems below.

\subsection{Comments on the proof technique}

The main lines of our argument follow \cite{delalande2022quantitative}, where the stability inequality for potentials is first established for the entropic semi-discrete optimal transport problem, and then extended to the Wasserstein distance by an approximation argument. However, the proofs in  \cite{delalande2022quantitative}  rely heavily on the fact that the cost is quadratic, so that up to a shift by removing the second moments, it is linear, hence both concave and convex. 

 Our first key realization is that the convexity assumption can actually be circumvented, and only concavity of the cost plays an essential role in the stability inequality for the Kantorovich potentials \cref{eq:stability-potential}. This leads to a more straightforward argument for the case $p=2$. The extension to general costs however requires further ideas. A crucial one is a generalization of the Pr\'ekopa-Leindler inequality to $\log$-concave measures, which is proved using a modification of McCann's optimal transport argument \cite{McCann}.  Actually, for $p \ge 2$, also the argument by Brascamp and Lieb \cite{brascamp1976extensions} can be used to obtain Prékopa-Leindler's inequality, however for $p<2$ McCann's argument seems necessary, since concavity of the cost fails (even allowing for a quadratic shift).  Curiously, our argument yields therefore a quantitative stability in optimal transport through an optimal transport derived inequality. A large part of the article, covering the $p\geq 2$ cases, is thus mostly a rewriting of the proof in \cite{delalande2022quantitative} : targeting the core of the arguments allows us to reduce the assumptions and therefore extend the results. The $1<p<2$ case follows the same line of proof but relies heavily on our generalization of Prékopa-Leindler's inequality.

 Finally, to obtain stability for the  optimal transport maps \cref{eq:stability-map}, a Gagliardo-Nirenberg type inequality is used in \cite{delalande2022quantitative}. For $p \ge 2$ the argument goes through similarly, however for $p<2$ this requires some modifications, in particular the use of fractional Sobolev spaces.

\subsection{Extensions and open problems}

We mention some extensions of our results that may be worth exploring, possibly using the techniques we developed.
\begin{enumerate}
\item
One notices immediately that \cref{thm:main-map} holds for target measures with bounded support, while \cref{thm:delalande} requires only bounds on their moments. It seems that this gap can be filled simply following the arguments from \cite{delalande2022quantitative}.

\item It may also be worth exploring whether quantitative stability results for potentials hold when the reference measure has unbounded support.  An interesting example may be that of a Gaussian density $\rho$. Let us remark however that extending the stability result for optimal transport maps to general $\log$-concave measures appears to be a more open-ended problem, as it would require establishing  suitable ``weighted'' Gagliardo-Nirenberg type inequalities. Some positive result in this direction has recently been obtained in \cite{letrouit2025gluing}.

\item We conjecture that quantitative stability bounds on Riemannian manifolds or even more general settings could be obtained using suitable versions of our arguments, at least under bounded curvature assumptions -- see e.g.\ \cite{cordero2006prekopa} for an extension of Prékopa–Leindler type inequalities on Riemannian manifolds. Stability results for optimal transport on Riemannian manifolds are the object of the recent article \cite{kitagawa2025stabriemann}.

\item In the case $p=2$, one can obtain stability in the optimization problem, namely $\Var_\rho( \phi - \phi_\mu) \lesssim \WW_\ell^2(\rho, \mu) -  \Big(\int \phi d\rho + \int \phi^\ell d\mu \Big)$, where $\phi_\mu$ denotes an optimal potential between $\rho$ and $\mu$ for the linear cost $\ell$ (see below for precise definitions). It may be worth exploring what happens when $p\neq2$.

\item For applications, it may be interesting to inquire in more detail the sharpness of exponents and the optimality of constants in \cref{thm:main-pot} and \cref{thm:main-map}.

\end{enumerate}

\subsection{Structure of the paper}

\cref{sec:DN} contains definitions and general well-known results. \cref{sec:pot} is devoted to the proof of \cref{thm:main-pot}, while \cref{sec:maps} focuses on the derivation of \cref{thm:main-map}. In both sections, we naturally split the arguments in three cases $p=2$, which we believe is particularly useful to understand the key points in our argument, $p >2$ which is a relatively easy variation of the former, and then the hardest case $1<p<2$, which requires substantially new tools.

\section{Notation and basic facts}\label{sec:DN}

\subsection{General notation}
Throughout the paper, we let $d\ge1$ and write $|\cdot|$ for the Euclidean norm on $\R^d$ and $\ang{\cdot, \cdot}$ for the scalar product.  Given a vector $v \in \R^d$ and $\alpha>0$, we write $v^{(\alpha)}:= |v|^{\alpha-1} v$. We let $\XX$, $\YY \subseteq \R^d$ denote compact sets and write $\CC(\XX)$ for the set of continuous functions on $\XX$, $\PP(\XX)$ for the set of Borel probability measures on $\XX$, $R_\XX$ the smallest positive real number such that $\XX \subseteq B(0,R_\XX)$, and $\diam(\XX)$ for the diameter of $\XX$. Unless otherwise stated, we always tacitly assume that $\rho \in \PP(\XX)$ and $\mu, \nu \in \PP(\YY)$.  Given $\psi \in \CC(\YY)$ and $\mu \in \PP(\YY)$, we write $\ang{\mu|\psi} = \int_{\YY} \psi d \mu$ for the duality pairing (which we also use for signed measures). When $\YY$ is finite, $\YY = \cur{y_i}_{i=1}^n$, we naturally identify functions $\psi \in \CC(\YY)$ with vectors $(\psi_i)_{i=1}^n = (\psi(y_i))_{i=1}^n \in \R^n$ and probabilities $\PP(\YY)$ with the simplex, i.e., $(\mu_i)_{i=1}^d = \cur{\mu(y_i)}_{i=1}^n \in [0,1]^n$ with $\sum_{i=1}^n \mu_i = 1$. We write  $\nabla = (\partial_i)_{i=1}^n$ for  the gradient operator on $\CC(\YY)$ induced by such identification with $\R^n$, and similarly $\nabla^2 = (\partial_{i,j}^2)_{i,j=1}^n$ for the Hessian. 

We always tacitly identify any absolutely continuous measure (with respect to the Lebesgue measure) with its density. In particular, we say that $\lambda$ is $\log$-concave if it has a density $\lambda = e^{-V}$ for some convex lower-semicontinuous $V: \R \to (- \infty, \infty]$. The support of a $\log$-concave measure is necessarily convex. Our main result for potentials may apply as well to singular $\log$-concave measures, arguing by approximations via convolutions : it in fact provides unicity of the entropic selection of an optimal potential (up to constants) but we leave this discussion to further investigations.

We write $\nor{\cdot}_{L^p(\rho)}$ for the $L^p$ norm with respect to the measure $\rho$ and simply $\nor{\cdot}_p$ if $\rho$ is understood. We occasionally use the notation $\nor{\cdot}_\infty$ for the uniform norm. Similarly, we write 
\begin{equation}
\mean_{\rho}(f) = \int_\XX f d \rho,
\end{equation}
\begin{equation}
 \Var_\rho(f) = \int_{\XX} \bra{ f - \int_{\XX} f d \rho}^2  d \rho
\end{equation}
for the mean and  variance of $f: \XX \to \R$ -- assuming that the $\rho$-mean of $f$ is well-defined. We recall the characterization
\begin{equation}
 \Var_\rho(f) = \min_{\lambda \in \R} \int_{\XX} \bra{ f - \lambda}^2  d \rho,
\end{equation}
which easily yields that if $\tilde \rho \in \PP(\XX)$ is absolutely continuous with respect to $\rho$, with a uniformly bounded density $d \tilde \rho /d \rho$, then
\begin{equation}\label{eq:variance-change-of-measure}
 \Var_{\tilde \rho} (f) \le \bra{ \sup_{\XX} \frac{d \tilde \rho}{d \rho}}  \Var_{\rho} (f).
\end{equation}
Moreover, $f\mapsto \Var_\rho(f)$ is easily seen to be convex. We write
\begin{equation}
  \osc(f) := \sup_{x \in \XX} f(x) - \inf_{x \in \XX} f(x)
\end{equation}
for the oscillation of $f$. We notice the elementary inequality
\begin{equation}
 \osc(f) \le \diam(\XX) \Lip(f),
\end{equation}
where $\Lip(f)$ denotes the Lipschitz constant of $f$ (possibly $\infty$ if $f$ is not Lipschitz continuous).

\subsection{Optimal transport and its entropic regularization}

Given a continuous cost function $c: \XX \times \YY \to \R$ 
we define the optimal transport problem with cost $c$ as the following minimization problem:
\begin{equation}\tag{$\text{OT}_c$}\label{pb:ot}
\TT^c(\rho, \mu) :=  \inf\limits_{\pi \in \Pi(\rho, \mu)} \ds\int\limits_{\XX \times \YY} c(x,y) d\pi(x,y)
\end{equation}
where $\Pi(\rho, \mu)$ 
denotes the set of couplings between $\rho$ and $\mu$. For any $\eps > 0$, we introduce the $\eps$-entropic regularized optimal transport problem, with cost $c$, as the following minimization problem:
\begin{equation}\tag{$\text{EOT}_{c,\eps}$}\label{pb:eot}
\TT^{c,\eps}(\rho, \mu) := \inf\limits_{\pi \in \Pi(\rho, \mu)}  \ds\int\limits_{\XX \times \YY}c(x,y) d\pi(x,y) + \eps \ent(\pi || \rho \otimes \mu)
\end{equation}
where $\ent$ denotes the relative entropy functional.

Note also that for $\sigma \in \PP(\YY)$ such that $\mu \ll \sigma$ and $\sigma \ll \mu$, it admits the following rewriting 
\begin{equation}
\TT^{c,\eps}(\rho, \mu) := \inf\limits_{\pi \in \Pi(\rho, \mu)}  \ds\int\limits_{\XX \times \YY}c(x,y) d\pi(x,y) + \eps \ent(\pi || \rho \otimes \sigma) - \eps \ent(\mu || \sigma),
\end{equation}
which will be quite useful in our derivations. Although $\sigma \ll \mu$ is not directly needed for the above equality to hold, it allows to define optimal potentials on the whole $\text{supp}(\sigma)$ which will prove useful in the following. From now on, $\sigma$ will always define an element of $\PP(\YY)$ satisfying $\mu \ll \sigma$ and $\sigma \ll \mu$.

Since for simplicity we assume that $\XX$, $\YY$ are compact and $c$ is continuous, it always holds $\TT^c(\rho, \mu) < \infty$ for every $\rho \in \PP(\XX)$, $\mu \in \PP(\YY)$. When $p \geq 1$ and  $c$ is the so-called $p$-cost,  $c(x,y) = \frac{1}{p} |x-y|^p$ for $x \in \XX$, $y \in \YY$, we write $\TT^p$ instead of $\TT^c$ and notice that
\begin{equation}\label{pb:otp}
 \TT^p(\rho, \mu) = \frac 1 p \WW_p(\rho,\mu)^p.
\end{equation}
When $p=2$,  writing $|x-y|^2 = |x|^2 + |y|^2 - 2 \ang{x,y}$, we collect the identity
\begin{equation}\label{eq:quadratic}
 \TT^2(\rho, \mu) = \frac 1 2  \int_{\XX} |x|^2 d\rho(x)+ \frac 1 2\int_{\YY} |y|^2 d \mu(y) + \TT^\ell (\rho, \mu),
\end{equation}
where $\ell$ denotes the linear cost $(x,y) \mapsto - \ang{x,y}$.

We will make crucial use of the following result on the so-called displacement interpolation between two absolutely continuous probabilities, which is well-known for the quadratic case (see e.g.\ McCann's proof of Prékopa-Leindler \cite[Appendix D]{McCann}), while for general power cost a proof can be found as a byproduct of the proof of \cite[Proposition 9.3.9]{AGS} applied to the (opposite of) the differential entropy functional $\FF(\rho) = \int \rho \log (\rho)$.  

\begin{lem}\label{lem:boundgeo} Let $p>1$, $\rho_0$, $\rho_1 \in \PP(\XX)$ be absolutely continuous and let $T$ be the 
optimal transport map for $\WW_p(\rho_0, \rho_1)$. For $t \in [0,1]$, let $T_t(x) = (1-t) x+ tT(x)$ and set $\rho_t = (T_t)_\sharp \rho_0$ the displacement interpolation between $\rho_0$ and $\rho_1$. Then $\rho_t$ is absolutely continuous and
\begin{equation}
\rho_t(T_t(x)) \leq \rho_0(x)^{1-t}  \rho_1(T(x))^t \quad \text{for $\rho_0$-a.e.\ $x \in \XX$.}
\end{equation}
\end{lem}

We also recall the following upper bound for the quadratic Wasserstein distance between two densities, in terms of the $L^2$ norm of their difference, with respect to a $\log$-concave reference measure $\rho$.

\begin{lem}\label{lem:peyre}
 Let $\rho \in \PP(\R^d)$ be $\log$-concave, with finite Poincaré constant. Then, there exists $C = C(\rho)<\infty$ such that the following holds. For every $\mu_0 = f_0\rho$, $\mu_1 = f_1 \rho \in \PP(\R^d)$, both absolutely continuous with respect to $\rho$,
 \begin{equation}
  \WW_2(\mu_0, \mu_1)^2 \le \frac{C}{\operatorname{ess-inf} f_1} \nor{ f_1 - f_0}_{L^2(\rho)}^2.
 \end{equation}
\end{lem}
Inequalities of this kind  have been used by many authors, see e.g. \cite[Corollary 3]{peyre2018comparison}, \cite[Proposition 2.3]{ambrosio2019pde}, \cite[Theorem 2]{Le17} and can be proved combining the Benamou-Brenier formula for $\WW_2$ and the Poincaré-Wirtinger inequality with respect to the measure~$\rho$.

\begin{proof}[Proof of \cref{lem:peyre}] Starting from \cite[Corollary 3]{peyre2018comparison}, with the setting of the lemma, it holds 
 \begin{equation}
  \WW_2(\mu_0, \mu_1)^2 \le \frac{C}{\operatorname{ess-inf} f_1} \nor{ f_1 - f_0}_{\dot{H}^{-1}(\rho)}^2.
 \end{equation}
The lemma follows from an application of Poincaré inequality with respect to the measure $\rho$ in order to control the $\dot{H}^{-1}$ norm by the $L^2$ norm. 
\end{proof}

\subsection{Dual formulations and Kantorovich functional}

Since the ``source'' measure $\rho \in \PP(\XX)$ plays a distinguished role in our results, we write the dual formulation of the optimal transport problem in the following asymmetric way (also called semi-dual in \cite{delalande2022quantitative}):
 \begin{equation}
  \TT^c(\rho, \mu) = \sup_{\psi \in \CC(\YY)}  \cur{ \int_{\XX} \psi^c d \rho + \int_{\YY} \psi d \mu },
 \end{equation}
where
\begin{equation}
 \psi^c(x) := \inf_{y \in \YY} \cur{ c(x,y) - \psi(y)}
\end{equation}
denotes the $c$-transform of $\psi$. For the $\eps$-entropic regularized version, duality reads \cite{marino2020optimal, nutz2022entropic}
 \begin{equation}
  \TT^{c,\eps}(\rho, \mu) = \sup_{\psi \in \CC(\YY)}  \cur{ \int_{\XX} \psi^{c,\eps} d \rho + \int_{\YY} \psi d \mu } - \eps \ent(\mu || \sigma),
 \end{equation}
where the $(c,\eps)$-transform of $\psi$ is given by
\begin{equation}
 \psi^{c,\eps}(x) := - \eps \log\bra{ \int_{\YY} \exp\bra{ - \frac{ c(x,y) - \psi(y)}{\eps} } d \sigma}.
\end{equation}
We notice that this no longer depends on $\mu$. Being able to use the same definition of $(c,\eps)$-transform for different target measure $\mu$ is a key element of the proof.  Symmetrically, $c$-transform and $(c,\eps)$-transform can be defined for continuous functions on $\XX$. When a function in $\CC(\YY)$ can be written as a $c$-transform (resp. $(c,\eps)$-transform) we say it is $c$-concave (resp. $(c,\eps)$-concave). In this formulation, a $c$-concave or $(c,\eps)$-concave optimizer $\psi \in \CC(\YY)$ is called a Kantorovich potential (for $\TT^c(\rho, \mu)$ or $\TT^{c,\eps}(\rho, \mu)$),  and $(\phi, \psi) = (\psi^c, \psi) \in \CC(\XX) \times \CC(\YY)$ is called a pair of Kantorovich potentials (and similarly for $\eps > 0$).  We notice that such transforms define monotone decreasing maps, i.e., for $\psi, \tilde \psi \in \CC(\YY)$
\begin{equation}
 \psi(y) \le \tilde \psi(y) \quad \text{for every $y \in \YY$} \quad \Rightarrow \quad \psi^{c,\eps}(x) \ge \tilde \psi^{c,\eps}(x) \quad \text{for every $x \in \XX$,}
\end{equation} and are concave,
\begin{equation}
 ((1-t) \tilde \psi + t \psi)^{c,\eps}(x) \ge (1-t) \tilde \psi^{c,\eps}(x) + t \psi^{c,\eps}(x) \quad \text{for every $t \in [0,1]$, $x \in \XX$,}
\end{equation}
as a consequence of H\"older's inequality. Moreover, $(\psi+\lambda)^{c,\eps} = \psi^{c,\eps} - \lambda$ for every constant $\lambda \in \R$.

We introduce, for $\phi_0$, $\phi_1 \in \CC(\XX)$, the quantity
\begin{equation}
M_{\phi_0,\phi_1}^c =  \sup_{t \in [0,1]}\cur{\osc(\phi_t)},
\end{equation}
where $\phi_t := ((1-t)\phi_0^c + t\phi_1^c)^c$.

\begin{ex}
If the cost $c$ is Lipschitz continuous, then any function $\phi = \psi^c$ that is a $c$-transform is Lipschitz continuous.  In the case of the $p$-cost, we obtain $\Lip(\phi) \le (R_\XX + R_\YY)^{p-1}$, hence $\osc(\phi) \le 2R_\XX (R_\XX + R_\YY)^{p-1}$. In particular, we can always bound from above
\begin{equation}\label{eq:M-p-cost}
 M_{\phi_0, \phi_1}^c \le 2R_\XX (R_\XX + R_\YY)^{p-1}.
\end{equation}
\end{ex}

As in \cite{delalande2022quantitative} we introduce the $c$-Kantorovich functional
\begin{equation}
\psi \in \CC(\YY)\quad  \mapsto  \quad
\KK^{c}(\psi) = -\int_{\XX}  \psi^{c}(x) d\rho(x)
\end{equation}
and similarly the $(c,\eps)$-Kantorovich functional
\begin{equation}
\psi \in \CC(\YY)\quad  \mapsto  \quad
\KK^{c,\eps}(\psi) = -\int_{\XX}  \psi^{c,\eps}(x) d\rho(x).
\end{equation}
We remark that the notation does not highlight the dependence upon $\rho$, which will be clear from the context. Both functionals $\KK^c$, $\KK^{c,\eps}$ are convex, and the duality formulas read
\begin{equation}\label{eq:dualformclassicandentropic}\begin{split}
 \TT^c(\rho, \mu) & = \sup_{\psi \in \CC(\YY)}  \cur{ \ang{\mu|\psi} -\KK^{c}(\psi) }\\
  \TT^{c,\eps}(\rho, \mu) & = \sup_{\psi \in \CC(\YY)}  \cur{ \ang{\mu|\psi} -\KK^{c,\eps}(\psi) }- \eps \ent(\mu || \sigma).
  \end{split}
\end{equation}

\begin{rem}\label{rem:vanishing-entropy}
 With the above definitions, we have convergence of the regularized problems to the usual transport problem as $\eps\downarrow 0$, as discussed in \cite{nutz2022entropic}. 
 More precisely, since both $\XX$ and $\YY$ are bounded and the cost function is (uniformly) continuous, any family of functions that are $(c,\eps)$-transforms has a uniform modulus of continuity, also uniformly with respect to $\eps>0$. If we require that optimal potential $\phi = \psi^{c,\eps}$ have e.g.\ zero mean with respect to $\rho$, we get by Arzelà-Ascoli that the families are relatively compact, hence up to a subsequence they converge uniformly. Reversing the roles of $\phi$ and $\psi$ leads to compactness also for the optimizers $\psi$'s. If moreover uniqueness of the (zero mean) optimal potentials  holds, as is the case of $p$-costs for $p>1$ in our setting, one obtains uniform convergence of the optimizers (and also of the costs), otherwise one can still extract a uniformly converging subsequence.
 \end{rem}

\subsection{Semi-discrete problem}

It will be convenient to further approximate by assuming that $\YY$ is a finite set so that we can apply finite-dimensional calculus. The following result extends \cite[Lemma 2.7]{delalande2022quantitative} and allows us to pass from such semi-discrete case to the general case of $\YY$ compact.

\begin{lem}\label{lem:approx} Let $\XX \subset \R^d$ be compact and  convex and $\YY \subset \R^d$ be compact, $\rho \in \PP(\XX)$ be absolutely continuous, $\mu_0, \mu_1 \in \PP(\YY)$ and $c: \XX\times \YY \to \R$ be continuous. For $k \in \{0,1\}$,  assume that there exists a unique Kantorovich potential with zero $\rho$-mean  for $\mathcal{T}^c(\rho, \mu_k)$ $\phi_k \in \CC(\XX)$ and let $\psi_k$ its $c$-transform. Then, there exist sequences of finite sets $(\YY^n)_{n=1}^\infty$, discrete measures $(\mu_{k}^n)_{n=1}^\infty$, Kantorovich potentials $(\phi_{k}^n, \psi_k^n)_{n=1}^\infty$ 
such that, for $k\in \cur{0,1}$,
\begin{enumerate}
\item $\operatorname{supp}(\mu_{k}^n) = \YY^n$ for every $n \ge 1$,
\item  $\nor{\phi_{k}^n - \phi_k}_{\infty} \to 0$ and $\nor{\psi_{k}^n - \psi_k}_{\infty} \to 0$ as $n \to \infty$
\item   it holds $\langle \mu_{k}^n | \psi_{1}^n - \psi_{0}^n \rangle  \to  \langle \mu_k | \psi_1 - \psi_0 \rangle$ as $n \to \infty$.
\end{enumerate}
\end{lem}

Without the uniqueness assumption, one may still extract a converging subsequence to a pair of Kantorovich potentials.

\begin{proof} For any $n \geq 1$,  consider a finite partition $\YY = \underset{1 \leq i \leq n}\bigsqcup \YY_{i}^n$ with $\eps_n := \max\limits_{1 \leq i \leq n} \diam(\YY_{i}^n)$ infinitesimal as $n \to \infty$. Define $\YY^n = \cur{y_i^n}_{i=1}^n$, with  $y_i^n \in \YY_{i}^n$ and set
\begin{equation}
 \mu_{k}^n = \sum_{i=1}^n \sqa{  \bra{ 1- \frac{1}{n}} \mu_k(\YY_{i}^n) + \frac{1}{n^2} } \delta_{y_{i}^n},
\end{equation}
so that (1) is verified.

Letting $\phi_{k}^n$ be the Kantorovich potential for $\mathcal{T}^c(\rho, \mu_k^n)$ with zero $\rho$-mean, then (2) holds by Arzelà-Ascoli and the fact that the limiting potential is unique.

To show (3), notice that 
\begin{equation}
  \WW_1(\mu_k, \mu_{k}^n) 
  \leq \eps_n + \frac{\diam(\YY)}{n},
\end{equation}
so that $\mu_{k}^n \to \mu_k$ weakly 
and therefore the limit holds, since the functions $\psi^n_k$ converge uniformly. 
\end{proof}

In the case where $\YY$ is finite, and fixing $\sigma$ as the uniform measure on $\YY$, we recall some formulas for the derivatives of $\KK^{c,\eps}$ which is easily seen to be a smooth map from $\CC(\YY)$ (identified with $\R^n$) to $\R$. These are proved in \cite[Lemma 4.2, Proposition 3.6]{delalande2022quantitative} for the case of the ``linear'' cost $\ell$ defined above \cref{eq:quadratic}, but the proof can straightforwardly be generalized to any cost. 
\begin{equation}
\partial_i (\psi^{c,\eps})(x)  = - \ds\frac{\exp\bra{ \frac{ \psi_i -c(x,y_i) }{\eps}}\sigma_i}{\sum\limits_{j=1}^n \exp\bra{\frac{\psi_j - c(x,y_j)}{\eps}}\sigma_j}, \quad \text{for $i=1\ldots, n$,}
\end{equation}
which we notice defines for every $x \in \XX$ a probability vector, i.e., $-\partial_i (\psi^{c,\eps})(x)  \in [0,1]$ and $-\sum_i \partial_i (\psi^{c,\eps})(x)  =1$, that we naturally identify with a probability on $\YY$. We can thus let 
\begin{equation}
\pi^c_{\eps,\psi}(y_i|x) := - \partial_i (\psi^{c,\eps})(x),
\end{equation}
\begin{equation}
\pi^c_{\eps,\psi}(x,y) := \pi^c_{\eps,\psi}(y|x)\rho(x)
\end{equation}
and finally 
\begin{equation}
\mu(y) = \int_{\XX} \pi^c_{\eps,\psi}(x,y) dx.
\end{equation}
Note that $\pi^c_{\eps,\psi}(\cdot | \cdot)$ is the decomposition against $\rho$ of $\pi^c_{\eps,\psi}$ and $\mu$ the second marginal of $\pi^c_{\eps,\psi}$. One can check that $\pi^c_{\eps,\psi}$ is the optimal entropic transport plan between $\rho$ and $\mu$ and $(\psi^{c,\eps}, \psi)$ a pair of optimal entropic potentials. Note that $\pi^c_{\eps,\psi}(\cdot | \cdot)$ will in turn

play a role in the expression for the gradient and  Hessian of $\KK^{c, \eps}$: 
\begin{equation}\label{eq:hessian-K}
\begin{split}
\ang{ \nabla \KK^{c,\eps}(\psi),  v} = \int_{\XX}  \mean_{\pi^c_{\eps,\psi}(.|x)} (v) d \rho(x), \\
 \langle v,  \nabla^2 \KK^{c,\eps}(\psi) v \rangle =  \frac{1}{\eps}\int_\XX   \Var_{\pi^c_{\eps,\psi}(.|x)}(v)  d\rho(x).
\end{split}
\end{equation}
In particular, when evaluated in $\psi^{c,\eps}$ the optimal potential for the transport problem between $\rho$ and $\mu$, it holds
\begin{equation}
\nabla \KK^{c,\eps}(\psi^{c,\eps}) = \mu.
\end{equation}
This can also be seen directly from the fact that $\psi^{c,\eps}$ is a critical point of the functional maximized on the second line of \cref{eq:dualformclassicandentropic}. 
We also introduce the following ``partition function''
\begin{equation}\label{eq:Z}
\psi \in \CC(\YY) \mapsto \I_\beta (\psi) := \int_\XX e^{\beta \psi^{c, \eps}} d\rho.
\end{equation}
As with the function $\KK$, the notation does not highlight the dependence upon $\rho$, which will be clear from the context. The function $\I$ is smooth and in particular the following formula
\begin{equation}\label{eq:hessI}
\langle v, \nabla^2 \log \I_\beta(\psi) v \rangle = -  \frac{\beta}{\eps} \int_{\XX} \Var_{\pi^c_{\eps,\psi}(.|x)} (v)   d \rho_{\beta, \psi}^\eps(x)  + \beta^2 \Var_{\rho_{\beta, \psi}^\eps}\bra{ \mean_{\pi^c_{\eps,\psi}(.|x)}(v) },
\end{equation}
holds with
\begin{equation}\label{eq:defrhoI}
\rho_{\beta,\psi}^\eps := \frac{1}{\I_\beta(\psi)} e^{\beta\psi^{c ,\eps}} \rho. 
\end{equation}

\begin{rem}
In most proofs, to keep notation simple we will fix $\beta = 1$ and write $\ZZ := \ZZ_1$, $\rho_{\psi}^\eps := \rho_{1,\psi}^\eps$.
\end{rem}

Finally, recalling that we denote with $\ell$ the cost $(x,y) \mapsto - \ang{x,y}$, we obtain
\begin{equation}
 \psi^\ell (x) = - \sup_{y\in \YY} \cur{\ang{x,y} + \psi(y)}, \quad  \psi^{\ell,\eps} (x) = -\eps \log \bra{ \int_{\YY} \exp\bra{ \frac{\ang{x,y} + \psi(y)}{\eps} }  d \sigma(y)}.
\end{equation}

\begin{rem}\label{prop:linkotctilde} Let $f \in L^1(\rho)$ and ${\tilde c}(x,y) = c(x,y) + f(x)$. Then the following equalities hold, with equality of the set of minimizers.
\begin{equation}
\TT^c(\rho, \mu) = \TT^{\tilde{c}}(\rho, \mu) + \mean_{\rho}(f) 
\end{equation}
\begin{equation}
\TT^{c,\eps}(\rho, \mu) = \TT^{\tilde{c},\eps}(\rho, \mu) + \mean_{\rho}(f) 
\end{equation}
Moreover, for any $\psi \in \CC(\YY)$,
\begin{equation}
\KK^{c,\eps}(\psi) = \KK^{\tilde{c},\eps}(\psi) + \mean_{\rho}(f) 
\end{equation}
In particular, if $(\phi_c, \psi_c)$ (resp. $(\phi_{c,\eps}, \psi_{c,\eps})$) denotes a couple of optimal potentials for the cost $c$ for \cref{pb:ot}, (resp. for \cref{pb:eot}),
\begin{equation}
(\phi_{\tilde{c}}, \psi_{\tilde{c}}) = (\phi_c +  f, \psi_c),
\end{equation}
\begin{equation}
(\phi_{\tilde{c},\eps}, \psi_{\tilde{c},\eps}) = (\phi_{c,\eps} + f, \psi_{c,\eps}).
\end{equation}
\end{rem}

\subsection{Fractional Sobolev spaces}

We start by recalling a definition of the spaces $W^{s,p}(\Omega)$, when the order $s$ is not an integer and $\Omega \subseteq \R^d$  is a bounded domain with Lipschitz boundary. For $0 < s < 1$ and $1 \leq p < \infty$, we set
\begin{equation}
 [u]_{s,p}^p :=  \ds\int_\Omega \int_\Omega \frac{ |u(x) - u(y)|^p}{|x-y|^{d + sp}} dxdy 
\end{equation}
and 
\begin{equation}
\|u\|_{W^{s,p}} := \|u\|_{L^p(\Omega)} +  [u]_{s,p} \qquad W^{s,p}(\Omega) := \{ u \in L^p(\Omega) \text{ and } \|u\|_{W^{s,p}} < \infty \}.
\end{equation}

When $s>1$ is not an integer, we write $s = m + \sigma$, with $m \in \N$, and $\sigma \in (0,1)$ and let
\begin{equation}
W^{s,p}(\Omega) := \{ u \in W^{m,p}(\Omega); D^mu \in W^{\sigma,p}(\Omega) \}
\end{equation}
normed with 
\begin{equation}
\|u\|_{W^{s,p}} := \|u\|_{W^{m,p}} + [D^mu]_{W^{\sigma,p}}.
\end{equation}
When $p=2$ we usually note $H^s(\Omega) := W^{s,2}(\Omega)$.

We will need the following lemma in the sequel. It is merely an iterated application of Gagliardo-Nirenberg inequality and even though we will only use it for $\Omega$ open interval of $\R$, we state it in full generality.
\begin{lem} For any $r \in (1,2)$ and any $u \in L^2(\Omega) \cap W^{1,\infty}(\Omega) \cap W^{r,1}(\Omega)$, $u \in H^1(\Omega)$ and the following inequality holds :
\begin{equation}
 \label{eq:gagliardo-useful}
 \|u\|_{H^1(\Omega)} \le C \|u\|_{L^2(\Omega)}^{1- \frac{2}{1+r}} \|u\|_{W^{1,\infty}(\Omega)}^{\frac{1}{1+r}} \nor{u}_{W^{r, 1}(\Omega)}^{\frac{1}{1+r}}.
\end{equation}
\end{lem}

\begin{proof}

We recall the fractional Gagliardo-Nirenberg inequality (see \cite[Theorem 1]{BreMir}):
for any $0 \le s_1<s_2 $, $p, p_1, p_2 \in [1, \infty]$, with $\theta \in (0,1)$,
\begin{equation}
s := \theta s_1 + (1-\theta)s_2, \quad \text{and} \quad \frac 1 p  = \theta \frac 1 {p_1} + (1-\theta) \frac 1 {p_2},
\end{equation}
then
\begin{equation}\label{eq:gn}
\|u\|_{W^{s,p}(\Omega)} \le C \|u\|_{W^{s_1,p_1}(\Omega)}^\theta  \|u\|_{W^{s_2,p_2}(\Omega)}^{1-\theta}
\end{equation}
provided that the following condition \emph{fails}:
\begin{equation}\label{eq:brezis-1}
 \text{$s_2$ is  an integer, $p_2=1$  and $s_1- \frac{1}{p_1} \ge s_2- \frac{1}{p_2}$.}
\end{equation}
By iterations, we may also consider $\theta_1$, $\theta_2$, $\theta_3 \in (0,1)$ with $\theta_1+\theta_2+\theta_3 = 1$, $0 \le s_1 < s_2 <s_3$ and $p, p_1, p_2, p_3 \in [1, \infty]$ such that
\begin{equation}\label{eq:brezis-2}
 s := \sum_{i=1}^3 \theta_i s_i, \quad \frac 1 p = \sum_{i=1}^3\frac{ \theta_i}{p_i},
\end{equation}
provided that both \cref{eq:brezis-1} and the following condition \emph{fail}:
 \begin{equation}
 \text{$s_3$ is an integer, $p_3=1$  and $ s- \frac{1}{p} \ge s_3- \frac{1}{p_3}$.}
\end{equation}
In our application, we need the following special case: $s_1 = 0 < s_2 = 1 < s_3<2$ and $p_1=2$, $p_2=\infty$, $p_3=1$, and $s=1$, $p=2$, so that $\theta_1$, $\theta_2$, $\theta_3$ must satisfy $\theta_1+\theta_2+\theta_3=1$ and
\begin{equation}
 1 = \theta_2 + \theta_3 s_3 \quad \text{and} \quad \frac 1 2 = \frac{\theta_1}{2}+ \theta_3.
\end{equation}
These conditions yield
\begin{equation}
 \theta_1 = 1- \frac{2}{1+s_3}, \quad \theta_2= \theta _3  = \frac{1}{1+s_3}.
\end{equation}
which must be also complemented with the failure of both \cref{eq:brezis-1}, which fails because $p_2=\infty$ and \cref{eq:brezis-2}, which fails because $s_3$ is not an integer. In conclusion, writing $r\in (1,2)$ instead of $s_3$, we have the desired inequality.
\end{proof}

\section{Quantitative stability of potentials}\label{sec:pot}

Aim of this section is to establish our first main result, \cref{thm:main-pot}.  To simplify our argument and highlight the differences with the results in \cite{delalande2022quantitative}, we assume that $\rho=\lambda$ is a $\log$-concave probability measure with bounded (convex) support and point out in \cref{rem:density} below the minor modifications for the general case of $\rho$ equivalent to $\lambda$. We also split our argument into three subsections. We deal first with the quadratic case, in particular recovering and slightly extending the results in \cite{delalande2022quantitative}, but mainly in order to point out where our arguments differ and allow for extensions respectively to the case $p\ge 2$ (\cref{sec:p>2})  and $1<p<2$ (\cref{sec:p<2}).

\subsection{The quadratic case}\label{sec:quadrcase}

In this section, we prove \cref{thm:main-pot} in the case $p=2$. In view of the identity \cref{eq:quadratic}, it is sufficient to argue in the case of the linear cost $c(x,y) = -\ang{x,y}$ that we denote with $\ell$. Therefore, we establish the following result.

\begin{theo}\label{theo:stabcompact}
Let $\rho$ be a $\log$-concave probability measure with bounded convex support $\XX\subseteq \R^d$ and let $\YY\subseteq \R^d$ be compact.  Given $\mu_0$, $\mu_1 \in \PP(\YY)$, let $(\phi_{\mu_0},\psi_{\mu_0})$, $(\phi_{\mu_1}, \psi_{\mu_1}) \in \CC(\XX) \times \CC(\YY)$ be pairs of Kantorovich potentials respectively for $\TT^\ell(\rho, \mu_0)$ and $\TT^\ell(\rho, \mu_1)$, then
\begin{equation}
\Var_\rho \bra{ \phi_{\mu_0} - \phi_{\mu_1}}  \leq 2M_{\phi_{\mu_0}, \phi_{\mu_1}}^\ell \langle \mu_0 - \mu_1 | \psi_{\mu_0}- \psi_{\mu_1} \rangle.
\end{equation}
\end{theo}

Recall that, by \cref{eq:M-p-cost} for $p=2$, we have $M_{\phi_{\mu_0}, \phi_{\mu_1}}^\ell \le 2R_\XX (R_\XX + R_\YY)$.

\begin{proof}[Proof of \cref{thm:main-pot} for $p=2$] Let us now prove \cref{thm:main-pot} from \cref{theo:stabcompact}. Being a Kantorovich potential for the $\ell$-cost, $\psi_{\mu_0}$ is the $\ell$-transform of $\phi_{\mu_0} \in \CC(\XX)$. As a consequence, $\psi_{\mu_0}$ is $R_\XX$-Lipschitz. Similarly, $\psi_{\mu_1}$ is $R_\XX$-Lipschitz. Then, by Kantorovich-Rubinstein duality formula, one gets
\begin{equation}
\la \mu_0 - \mu_1 | \psi_{\mu_0} - \psi_{\mu_1} \ra \leq 2R_\XX \WW_1(\mu_0,\mu_1)
\end{equation}
which, put along with \cref{theo:stabcompact} yields \cref{thm:main-pot}.

\end{proof}

%
%
%

We now turn to the proof of \cref{theo:stabcompact}. It is known or easy to check that while strict convexity yields uniqueness in a minimization problem, strong convexity yields explicit stability estimates. Following the idea developed in \cite{delalande2022quantitative}, we aim at a suitable convexity property for the functional $\KK^\ell$, which we obtain first in the $\eps$-entropic semi-discrete  (i.e., for finite $\YY$) case, then via an approximation argument for the general case of $\YY$ compact. However, to highlight the differences (and the improvements) with respect to \cite{delalande2022quantitative}, we split the proof of \cref{theo:stabcompact} into three steps:
\begin{enumerate}
 \item[\emph{Step 1 ($\log$-concavity of $\I$)}.] Arguing for fixed $\eps>0$ and $\YY$ finite, we establish $\log$-concavity of the ``partition function'' $\I$ introduced in \cref{eq:Z}.
Unlike \cite[Proposition 4.6]{delalande2022quantitative}, where this is obtained as a ``black box'' application of the Pr\'ekopa-Leindler inequality when $\rho$ is the Lebesgue measure on a convex bounded $\XX$, we obtain this property for a general $\log$-concave measure $\rho$ through an optimal transport argument (\cref{theo:Iconcav}), which itself can be used to prove the Pr\'ekopa-Leindler inequality, but is more amenable to further generalizations.

\item[\emph{Step 2 (modified convexity of $\KK^{\ell,\eps}$)}.] By direct computation, we  go from the $\log$-concavity of $\I$ to a strong convexity-like inequality on $\KK^{\ell, \eps}$, that reads for $\psi_0$, $\psi_1 \in \CC(\YY)$,
\begin{equation}
 \Var_{\rho}(\psi^{\ell,\eps}_1 - \psi^{\ell,\eps}_0)  \leq C \langle \nabla \KK^{\ell,\eps}(\psi_1) - \nabla \KK^{\ell, \eps}(\psi_0) | \psi_1 - \psi_0 \rangle,
\end{equation}
for a suitable  constant $C$. Let us stress that the inequality above provides  quantitative convexity, but not strong convexity in the usual sense, for the dual potentials appear in the left hand side. Still, our use of integration against $\rho$ in the first step allows us to yield a constant independent of upper and lower bounds on the density of $\rho$ with respect to the Lebesgue measure, going beyond the case covered in \cite{delalande2022quantitative}. Moreover, when compared to \cite{delalande2022quantitative}, where the (stronger) inequality with $\Var_{(\mu+\nu)/2}(\psi_\mu-\psi_\nu)$ in the left hand side is established, our proof does not rely on linearity of the cost (while \cite[Lemma 2.3]{delalande2022quantitative} apparently does).

\item[\emph{Step 3 (approximation and scaling)}] This final step consists in moving from the $\eps$-entropic regularization to the classical optimal transport by letting $\eps \downarrow 0$, as briefly discussed in \cref{rem:vanishing-entropy}, and then via an approximation argument to go from the discrete case of finite $\YY$ to the general  case of  compact $\YY$ using \cref{lem:approx}. Finally,  thanks to a scaling argument we improve the constant in the stability inequality: indeed, it is sufficient to introduce an ``inverse temperature'' parameter $\beta>0$, defining
\begin{equation}
 \I_\beta( \psi):=  \int_\XX e^{\beta \psi^{\ell, \eps}} d\rho
\end{equation}
and, after repeating the first two steps, we optimize upon $\beta$. This step does not differ substantially from \cite[Section 4.5, Section 2.4, Proposition 3.3]{delalande2022quantitative}. %
The scaling argument, also in \cite{delalande2022quantitative}, is actually of great importance for the case $1<p<2$ and will be detailed in  \cref{sec:p<2}, while adapting it to the present setting is straightforward.
\end{enumerate}

In view of the scheme detailed above, we establish the following proposition (settling \emph{Step 1}).
\begin{prop}[$\log$-concavity of $\I$]\label{theo:Iconcav} Let $\rho$ be a $\log$-concave measure with (convex) bounded support $\XX$, let $\YY\subseteq \R^d$ be finite, $\eps>0$ and $\beta>0$. For every $\psi_0$, $\psi_1 \in \CC(\YY)$ and $t \in [0,1]$, it holds
\begin{equation}\label{eq:logconcavI}
\log \I_\beta ((1-t)\psi_0 + t \psi_1) \geq (1-t) \log \I_\beta (\psi_0) + t \log \I_\beta (\psi_1).
\end{equation}
\end{prop}

\begin{proof} Recall that we identify $\CC(\YY)$ with $\R^n$, since $\YY = \cur{y_i}_{i=1}^n$, and write for brevity $\psi_i := \psi(y_i)$,  $\mu_i := \mu(y_i)$. To simplify the notation, we also simply write $\psi^\ell$ instead of $\psi^{\ell, \eps}$ and argue in the case $\beta=1$ only. Thus, given $\psi_0$, $\psi_1 \in \R^n$ and $x_0, x_1 \in \XX$ and setting $\psi_t = (1-t)\psi_0 + t\psi_1$ and $x_t = (1-t)x_0 + tx_1$, we have
\begin{equation}
 \begin{split}
 -(\psi_t)^\ell(x_t) &=  \eps \log\Big( \sum\limits_{i=1}^n \exp\bra{ \frac{\langle x_t, y_i \rangle + (\psi_t)_i}{\eps}} \sigma_i \Big)  \\
 &=  \eps \log\Big( \sum\limits_{i=1}^n \exp\bra{ \frac{(1-t) \langle x_0, y_i \rangle + t \langle x_1, y_i \rangle + (1-t) (\psi_0)_i + t (\psi_1)_i}{\eps}} \sigma_i \Big).
\end{split}
\end{equation}
Using H\"older's inequality, we get
\begin{equation}
- (\psi_t)^\ell (x_t) \leq  -(1-t) (\psi_0)^\ell(x_0)  - t  (\psi_1)^\ell(x_1)
 \end{equation}
  Therefore, taking the exponential and setting $h_t(x) = \exp((\psi_t)^\ell(x))$ we obtain the inequality
\begin{equation}\label{eq:hypPL}
h_t(x_t) \geq h_0(x_0)^{1-t} h_1(x_1)^t 
\end{equation}
As mentioned earlier, one can now use as a black box the Prékopa-Leindler inequality for $\log$-concave measures which yields $\|h_t\|_{L^1(\rho)} \geq  \|h_0\|_{L^1(\rho)} ^{1- t} \|h_1\|_{L^1(\rho)} ^t$. For completeness we recall the optimal transport version of the proof. Renormalizing, let us define $\hat{h}_k := h_k/ \|h_k\|_{L^1(\rho)}$ so that $\hat{h}_0$ and $\hat{h}_1$ are probability densities with respect to $\rho$.

By \cref{lem:boundgeo} applied with $p=2$, $\rho_0 := \hat{h}_0 \rho$, $\rho_1 := \hat{h}_1 \rho$, using also the $\log$-concavity of $\rho$ and writing $\rho_t$ the displacement interpolation between $\hat{h}_0\rho$ and $\hat{h}_1\rho$, $T$ the optimal map between $\hat{h}_0\rho$ and $\hat{h}_1\rho$, setting $T_t := (1-t) Id + tT$, we get for $\hat{h}_0\rho$-a.e.\ $x \in \XX$  the inequality
\begin{equation}
 \begin{split}
  \rho_t(T_t(x)) &\leq (\hat{h}_0(x)\rho(x))^{1-t} (\hat{h}_1(T(x))\rho(T(x)))^t \\
&\leq \hat{h}_0(x)^{1-t} (\hat{h}_1(T(x)))^t \rho(x)^{1-t}\rho(T(x))^t \\
&\leq \hat{h}_0(x)^{1-t} (\hat{h}_1(T(x)))^t \rho(T_t(x))\\
&\leq h_t(T_t(x))\rho(T_t(x)) \|h_0\|_{L^1(\rho)}^{-(1-t)} \|h_1\|_{L^1(\rho)}^{-t} .
 \end{split}
\end{equation}
Therefore, we deduce that
\begin{equation}
\|h_0\|_{L^1(\rho)}^{1-t} \|h_1\|_{L^1(\rho)}^{t}\rho_t(x) \leq  h_t(x)\rho(x) \quad \text{holds $\rho_t$ a.e.}
\end{equation}
Integrating we get the Prékopa-Leindler inequality :
\begin{equation}
\|h_t\|_{L^1(\rho)} \geq \|h_0\|_{L^1(\rho)}^{1-t}  \|h_1\|_{L^1(\rho)}^t.
\end{equation}
Taking the $\log$ and noticing that $\nor{h_t}_{L^1(\rho)} = \I(\psi_t)$  yields the thesis.

\end{proof}

In order to address \emph{Step 2} with lighter notations, fixing $\psi_0$, $\psi_1 \in \CC(\YY)$ and $\psi_t = (1-t)\psi_0 + t \psi_1$, we write $\pi_t(.|x) := \pi^\ell_{\eps,\psi_t}(.|x)$,
so that  \cref{eq:hessI} reads
\begin{equation}\label{eq:hessI-t}
\langle v, \nabla^2 \log \I(\psi_t) v \rangle = -  \frac{1}{\eps} \int_{\XX} \Var_{\pi_t(.|x)} (v)   d \rho_{\psi}^\eps(x)  + \Var_{\rho_{\psi}^\eps}\bra{ \mean_{\pi_t(.|x)}(v) },
\end{equation}

\begin{prop}\label{prop:Kstrongconv} Let $\rho$ be a $\log$-concave probability measure with compact support  $\XX \subset \R^d$, let $\mu \in \PP(\YY)$ with $\YY\subseteq \R^d$ finite and let $\eps > 0$. Given  $\psi_0$, $\psi_1 \in \CC(\YY)$, set
\begin{equation}
  \phi_t^\eps := ((1-t)\psi_0 + t\psi_1)^{\ell, \eps} \quad \text{for $t \in [0,1]$,}
\end{equation}

Then,
\begin{equation}\label{eq:notstrongconvex}
\Var_{\rho}(\phi_1^{\eps} - \phi_0^{\eps} )  \leq  C \langle \nabla \KK^{\ell,\eps}(\psi_1) - \nabla \KK^{\ell,\eps}(\psi_0) | \psi_1 - \psi_0 \rangle
\end{equation}
with $C := \exp(2 \sup_{t \in [0,1]} \osc(\phi_t^\eps))$.

\end{prop}

\begin{proof}
Again, we identify $\CC(\YY)$ with $\R^n$. From \cref{theo:Iconcav}, we get that $\nabla^2 \log \I(\psi_t) \in \R^{n\times n}$ is a symmetric negative semi-definite matrix, which yields for every $v \in \R^n$,
\begin{equation}
\ang{ v,  \nabla^2 \log \I(\psi_t) v }\leq 0  \quad  \forall \, 0 \leq t \leq 1.
\end{equation}
By the identity \cref{eq:hessI}, we get
\begin{equation}\label{eq:hessian-I-used}
 \Var_{\rho_{\psi_t}^\eps}\bra{ \mean_{\pi_t(.|x)}(v) } \le \frac{1}{\eps} \int_{\XX} \Var_{\pi_t(.|x)} (v)   d \rho_{\psi_t}^\eps  \quad  \forall \, 0 \leq t \leq 1.
\end{equation}
Our goal is now to change the integration against $\rho_{\psi_t}^\eps$ to integration against $\rho$.  Recalling the definition of $\rho_{\psi}^\eps$ in \cref{eq:defrhoI}, we deduce

\begin{equation}\label{ineq:boundrho}
e^{-\osc(\phi_t^\eps)} \rho_{\psi_t}^\eps \leq \rho \leq e^{ \osc(\phi_t^\eps)} \rho_{\psi_t}^\eps.
\end{equation}
We remark that here is where the choice of $\rho$ as a ``reference'' measure in the definition of $\I$, instead e.g., of the Lebesgue measure on $\XX$, turns out to be useful, as it frees us from making assumptions on the density of $\rho$ (besides its $\log$-concavity). Using \cref{eq:variance-change-of-measure} in both sides of \cref{eq:hessian-I-used} thus find, for every $0 \leq t \leq 1$,
\begin{equation}\label{eq:variance-bound}
 \Var_{\rho} \bra{\mean_{\pi_t(.|x)}(v)  } \leq \frac{  e^{2 \osc(\phi^\eps_t) }}{\eps} \int_\XX  \Var_{\pi_t(.|x)}  (v)d\rho = e^{2 \osc(\phi^\eps_t)}  \langle v,  \nabla^2 \KK^{\ell,\eps}(\psi_t) v \rangle
\end{equation}
where the second identity follows from \cref{eq:hessian-K}.

Recalling that $\phi_t^\eps := (\psi_t)^{\ell, \eps}$, we have by the chain rule
\begin{equation}\label{eq:derivphit}
\frac{d}{dt} \phi_t^{\eps}(x) = \ang{\nabla \psi_t^{\ell, \eps}(.|x) |  \psi_1 - \psi_0}  = - \mean_{\pi_t(.|x)}(\psi_1 - \psi_0)
\end{equation}
hence, by convexity of the variance,
\begin{equation}
\begin{split}
  \Var_{\rho}\bra{ \phi_1^{\eps} - \phi_0^{\eps} } & = \Var_{\rho}\bra{  \int_0^1\frac{d}{dt} \phi_t^{\eps} dt } \\
  & \le \int_0^1  \Var_{\rho} \bra{ \frac{d}{dt} \phi_t^{\eps}} dt \\
  & \le \int_0^1 \Var_{\rho} \bra{ \mean_{\pi_t(.|x)}(\psi_1 - \psi_0)} dt.
  \end{split}
\end{equation}
Using \cref{eq:variance-bound} with $v = \psi_1 - \psi_0$, we thus obtain
\begin{equation}
 \begin{split}
\Var_{\rho}(\phi_1^{\eps} - \phi_0^{\eps} )  &\leq \int_0^1  e^{2 \osc(\phi^\eps_t) } \langle \psi_1 - \psi_0,  \nabla^2  \KK^{\ell,\eps}(\psi_t) (\psi_1 - \psi_0) \rangle dt\\
& \leq C \int_0^1 \langle \psi_1 - \psi_0, \nabla^2  \KK^{\ell,\eps}(\psi_t) (\psi_1 - \psi_0) \rangle dt,\\
& \leq C \langle \psi_1 - \psi_0, \int_0^1 \nabla^2  \KK^{\ell,\eps}(\psi_t) (\psi_1 - \psi_0)dt \rangle
 \end{split}
\end{equation}
with $C= \exp(2 \sup_{t \in [0,1]} \osc(\phi_t^\eps))$. Finally, using that
\begin{equation}
 \nabla^2 \KK^{\ell,\eps}(\psi_t) (\psi_1-\psi_0) = \frac{d}{dt} \nabla \KK^{\ell,\eps}(\psi_t)
\end{equation}
we obtain the thesis.
\end{proof}
\begin{rem}\label{rem:sufficient_concav}
Let the reader note, for further use, that we actually only used
\begin{equation}
\la \psi_1 - \psi_0, \nabla \log\I(\psi_1) - \nabla \log \I( \psi_0) \ra = \int_0^1 \la \psi_1 - \psi_0, \nabla^2 \log\I(\psi_t), \psi_1 - \psi_0 \ra \leq 0.
\end{equation}
In particular, an inequality of the form 
\begin{equation}
\la \psi_1 - \psi_0, \nabla \log\I(\psi_1) - \nabla \log \I( \psi_0) \ra \leq C_{\psi_1, \psi_0}
\end{equation}
would similarly yield
\begin{equation}
\Var_{\rho}(\phi_1^{\eps} - \phi_0^{\eps} )  \leq  C \langle \nabla \KK^{\ell,\eps}(\psi_1) - \nabla \KK^{\ell,\eps}(\psi_0) | \psi_1 - \psi_0 \rangle + C^\frac12 C_{\psi_1, \psi_0},
\end{equation}
with $C$ as defined in \cref{prop:Kstrongconv}.
\end{rem}

 Let us recall here that we refer to  \cref{rem:vanishing-entropy} and \cref{lem:approx} for the approximation \emph{Step 3} and to  \cite[Section 4.5, Section 2.4, Proposition 3.3]{delalande2022quantitative} for detailed proofs. 

 \begin{rem}\label{rem:density}
 The only point where $\log$-concavity of $\rho$ is used is in the proof of \cref{theo:Iconcav}. Therefore, in order to cover the case of $\rho$ that admits a density uniformly bounded from above and below by strictly positive constants with respect to a $\log$-concave measure $\lambda$ with bounded support $\XX$,
 \begin{equation}
  m_\rho \lambda(x) \le  \rho(x) \le M_\rho \lambda(x) \quad \text{for $\lambda$-a.e.\ $x \in \XX$,}
 \end{equation}
 we need to define $\I$ using $\lambda$ instead of the measure $\rho$, namely $ \I_\beta( \psi):=  \int_\XX e^{\beta \psi^{\ell, \eps}} d\lambda $. This leads to an additional step in the argument, as $\log$-concavity of $\I$ leads to a modified version  of \cref{eq:hessian-I-used}, which reads
 \begin{equation}
 \Var_{\lambda_{\psi_t}^\eps}\bra{ \mean_{\pi_t(.|x)}(v) } \le \frac{1}{\eps} \int_{\XX} \Var_{\pi_t(.|x)} (v)   d \lambda_{\psi_t}^\eps  \quad  \forall \, 0 \leq t \leq 1.
\end{equation}
It is then sufficient to replace the variance and integration with respect to $\lambda_{\psi_t}^\eps$ with those with respect to $\rho_{\psi_t}^\eps$, so that the proof may then proceed without further changes: by the very definition of $\lambda_{\psi_t}^\eps$, we find the inequalities
\begin{equation}
  \frac{m_\rho}{M_{\rho}} \lambda_{\psi_t}^\eps (x) \le  \rho_{\psi_t}^\eps (x) \le  \frac{M_\rho}{m_{\rho}} \lambda_{\psi_t}^\eps (x) \quad \text{for $\lambda$-a.e.\ $x \in \XX$.}
\end{equation}
Therefore, recalling \cref{eq:variance-change-of-measure} we obtain that \cref{eq:hessian-I-used} holds with an additional factor $(M_\rho/m_\rho)^2$ in the right hand side. Similar arguments apply in the general cases we consider below, so we only focus on the case of a $\log$-concave measure $\rho$ to keep the exposition simple.
 \end{rem}

\subsection{The case $p\geq2$}\label{sec:p>2}

In this section, we prove \cref{thm:main-pot} in the case $p\geq2$. The previous proof does not fully use the linearity of the cost but only its concavity, hence an appropriate change of cost allows to generalize it with minimum effort.

\begin{theo}\label{theo:stabp>2}
Let $\rho$ be a $\log$-concave probability measure with bounded convex support $\XX\subseteq \R^d$ and let $\YY\subseteq \R^d$ be compact finally let $c: \XX \times \YY \to [0, + \infty)$ be a $C^2$ cost function. Given $\mu_0$, $\mu_1 \in \PP(\YY)$, let $(\phi_{\mu_0}, \psi_{\mu_0})$, $(\phi_{\mu_1},\psi_{\mu_1}) \in \CC(\XX) \times \CC(\YY)$ be pairs of Kantorovich potentials respectively for $\TT^c(\rho, \mu_0)$ and $\TT^c(\rho, \mu_1)$ then, there exists $C =C(c, \rho, \YY)<\infty$ such that
\begin{equation}
\Var_\rho \bra{ \phi_{\mu_0} - \phi_{\mu_1}}  \leq C \langle \mu_0 - \mu_1 | \psi_{\mu_0} - \psi_{\mu_1} \rangle.
\end{equation}
Moreover, for the $p$-cost one can take $C(p, \rho, \YY) = 4pR_\XX (R_\XX + R_\YY)^{p-1}$.
\end{theo}

\begin{proof}[Proof of \cref{thm:main-pot} for $p\ge2$] 
The thesis follows from \cref{theo:stabp>2} in the same way as \cref{thm:main-pot} for $p=2$ followed from \cref{theo:stabcompact}.
\end{proof}

%
%
%
%

This theorem follows mainly from the observation that concavity in the first variable is a sufficient condition on the cost for \emph{Step 1} to hold. \emph{Step 2} and  \emph{Step 3} will then apply as they hold for $c$ continuous and $\XX$ and $\YY$ compact. We conclude thanks to an appropriate modification of the cost. \emph{Step 1} is settled with the following proposition. It is equivalent to \cite[Lemma 5.3]{chizat2024sharper} and the proof follows the same lines.

\begin{prop}[$\log$-concavity of $\I$] \label{prop:csconcavI}  Let $\rho$ be a $\log$-concave measure and $\XX$ its support, let $\YY$ be finite. Then, the concavity of $x \to c(x,y)$ for all $y \in \YY$ is a sufficient condition for \cref{eq:logconcavI} to hold, for every $\psi_0$, $\psi_1 : \YY \to \R$.
\end{prop}

\begin{proof} Assuming that $\YY = \cur{y_i}_{i=1}^n$, and recalling that
\begin{equation} - (\psi_t)^{c,\epsilon}(x_t) = \epsilon \log\bra{ \sum\limits_{i=1}^n \exp\bra{\frac{(\psi_t)_i - c(x_t,y_i)}{\epsilon}}},
 \end{equation}
 with $x_t = (1-t)x_0 + tx_1$, one simply uses first the concavity of $c$ and then H\"older's inequality.
 \end{proof}

\begin{defi}\label{def:shiftc} For any $\CC^2$ cost function $c$, we introduce the shifted cost $\tilde{c}: \XX \times \YY \to \R_+$, defined by
\begin{equation}\label{eq:shiftc}
 \tilde{c}(x,y) = c(x,y) - \frac{\gamma}{2} |x|^2
\end{equation}
 where
\begin{equation} \gamma := \sup\limits_{x \in \XX,y \in \YY} \nor{ \nabla_x^2 c(x,y)}_{op},
\end{equation}
and $\nor{\cdot}_{op}$ denotes the operator norm.
\end{defi}
When $\XX$ and $\YY$ are two compact subsets of $\R^d$, the shifted cost is well defined and concave in the variable $x$.
In particular, this definition holds for the $p$-cost when $p \geq 2$, with $\gamma = (p-1)(R_\XX + R_\YY)^{p-2}$. Note also that this is actually the same exact strategy than going from the quadratic optimal transport problem to the maximum correlation problem where we go from a convex cost to a linear (and thus concave) one. In the current case, using the modified convexity \cref{eq:notstrongconvex} and not usual strong convexity is crucial, because the cost is no longer linear.

\begin{proof}[Proof of \cref{theo:stabp>2}] Letting $\XX$, $c$, $\rho$ as in the theorem and setting ${\tilde c}$ as in \cref{def:shiftc}, \cref{prop:linkotctilde} shows that proving \cref{theo:stabp>2} for ${\tilde c}$ is equivalent to proving it with $c$. We can therefore suppose $c$ concave. Then for $\YY$ finite \emph{Step 1} applies thanks to \cref{prop:csconcavI}. \emph{Step 2} as proved in \cref{prop:Kstrongconv} still holds. For \emph{Step 3}  we argue letting first $\eps \downarrow 0$ by \cref{rem:vanishing-entropy} and then use \cref{lem:approx} to extend from the semi-discrete to the general case. Note finally that the only requirement on $C$ is $C(p, \rho, \YY) \ge 2 M^{\tilde c}_{\phi_0, \phi_1}$, for all $\phi_0$, $\phi_1$ Kantorovich potentials and where $M^{\tilde c}$ means that it is computed with the ${\tilde c}$ cost. It is then easy to obtain a bound similar to \cref{eq:M-p-cost} and thus completing the thesis.
\end{proof}

\begin{rem}
 The entire argument straightforwardly extends to cost functions $c$ that are Lipschitz continuous and such that, for some $\gamma \in \R$, the modified cost \cref{eq:shiftc} is concave with respect to the variable $x$, for every $y \in \YY$. This larger family includes for example the ``boundary cost'' from \cite{figalli2010new},
\begin{equation}\label{eq:boundary-2}
b_{\Omega}^2(x,y) := \min( |x-y|^2, d(x,\Omega^c)^2 + d(y, \Omega^c)^2),
\end{equation}
where  $\Omega \subseteq \R^d$ and $d(\cdot, \Omega^c) = \inf_{z \in \Omega^c} |\cdot-z|$. 
\end{rem}

\subsection{The case $1<p<2$}\label{sec:p<2}

We finally move to the most challenging case of exponents $p \in (1,2)$ for which in particular, the smoothness or even the concavity condition of the above remark fail. We prove the following inequality.

\begin{theo}\label{theo:finalestp}
Let $p\in (1,2)$, set $q=p/(p-1)$, let $\rho$ be a $\log$-concave probability measure with bounded (convex) support $\XX\subseteq \R^d$ and let $\YY\subseteq \R^d$ be compact. Given $\mu_0$, $\mu_1 \in \PP(\YY)$, let $(\phi_{\mu_0}, \psi_{\mu_0})$, $(\phi_{\mu_1},\psi_{\mu_1}) \in \CC(\XX) \times \CC(\YY)$ be pairs of Kantorovich potentials respectively for $\WW_p(\rho, \mu_0)$ and $\WW_p(\rho, \mu_1)$. Then, there exists $C =C(p, \rho, \YY)<\infty$ such that
\begin{equation}
\Var_\rho(\phi_1 - \phi_0) \leq C \langle \mu_1 - \mu_0 | \psi_1 - \psi_0 \rangle^{\frac 2 q}.
\end{equation}
\end{theo}

As already noticed, we cannot restrict ourselves to the previous $\mathcal{C}^2$ case as the cost does not have bounded derivatives on the diagonal. We will derive a weaker version of \emph{Step 1}, before performing \emph{Step 2} and \emph{Step 3}, where  the scaling argument will be more important to control the additional term.

\begin{prop}\label{prop:Iconcavp}  Let $\rho$ be $\log$-concave with (convex) bounded support $\XX$, $\YY\subseteq \R^d$ be finite, $\eps>0$ and $\beta>0$. For every $\psi_0$, $\psi_1 \in \CC(\YY)$ and $t \in [0,1]$, it holds
\begin{equation}\label{eq:logconcavIp}\begin{split}
\log \I_\beta ((1-t)\psi_0 + t \psi_1) & \geq (1-t) \log \I_\beta (\psi^0) + t \log \I_\beta (\psi^1) \\
& \quad - \beta t(1-t)\gamma \WW_p(\rho_{\beta, \psi_0}^\epsilon,\rho_{\beta, \psi_1}^\epsilon)^p.
\end{split}
\end{equation}
 where $\gamma = \gamma(p)<\infty$ and for $k \in \cur{0,1}$, $\rho_{\beta, \psi_k}^\epsilon$ is defined as in \cref{eq:defrhoI} with $\psi_k$ instead of $\psi$. 
\end{prop}

Before the proof, let us collect some vectorial inequalities.

One can find the following in \cite[Chapter 12, (III)]{lindqvist2019notes}
\begin{equation}\label{eq:tooltocurvature}
p \langle |a -z|^{p-2}(a-z) - |b-z|^{p-2}(b-z) ,  a -b \rangle \leq \gamma |a - b|^p,
\end{equation}
valid for all $a,b,z \in \R^d$, since $1<p<2$ and for some constant $\gamma = \gamma(p) >0$. 

For the sake of completeness, let us derive this inequality. It is easily seen that the case $z=0$ is sufficient and that with an application of Cauchy-Schwarz inequality, it is enough to prove :
\begin{equation}
| |a|^{p-2} a - |b|^{p-2}b| \leq \gamma |a-b|^{p-1}.
\end{equation}
Starting from
\begin{equation}
| |a|^{p-2} a - |b|^{p-2}b|^2 = | |a|^{p-1} - |b|^{p-1}|^2 + 2|a|^{p-2}|b|^{p-2}(|a||b| - \la a,b \ra),
\end{equation}
we bound the first term, as $t \to t^{p-1}$ is $(p-1)$-Hölder continuous on $\R_+$ and by triangular inequality 
\begin{equation}
\begin{split}
| |a|^{p-1} - |b|^{p-1}|^2 & \leq \Big( \Big| |a| - |b| \Big|^{p-1} \Big)^2 \\
& \leq |a-b|^{2(p-1)}
\end{split}
\end{equation}
and the second, assuming $a,b \neq 0$, 
\begin{equation}
\begin{split}
|a|^{p-2}|b|^{p-2}(|a||b| - \la a,b \ra) & \leq |a|^{p-1}|b|^{p-1}(1 - \frac{\la a,b \ra}{|a||b|})^{p-1} (1 - \frac{\la a,b \ra}{|a||b|})^{2-p} \\
& \leq 2^{2-p} (|a||b| - \la a,b \ra)^{p-1} \\
& \leq 2^{3-2p} |a-b|^{2p-2}
\end{split}
\end{equation}
where we used for the last line $|a||b| - \la a,b \ra \leq \frac{|a-b|^2}{2}$. Both bounds together yields \cref{eq:tooltocurvature} in the case $z=0$.

We now deduce from it a semi-concavity inequality for $|.|^p$ when $1<p<2$. 
Convexity of $x \to |x|^p$ gives the two following inequalities, for $x_0, x_1 \in \R^n$, letting $x_t := (1-t)x_0 + t x_1$~: 
\begin{equation}
\begin{split}
|x_t|^p &\geq |x_0|^p + p \langle |x_0|^{p-2} x_0 | t(x_1 - x_0) \rangle \, \text{ and } \\
|x_t|^p &\geq |x_1|^p + p \langle |x_1|^{p-2} x_1 | (1-t) (x_0 - x_1) \rangle.
\end{split}
\end{equation}
Combining both, we get
\begin{equation}
|x_t|^p \geq (1-t) |x_0|^p + t|x_1|^p + pt(1-t) \langle |x_0|^{p-2}x_0 - |x_1|^{p-2} x_1 | x_1 -x_0 \rangle,
\end{equation}
and finally, applying \cref{eq:tooltocurvature} :
\begin{equation}\label{eq:pPC}
|x_t|^p \geq (1-t) |x_0|^p + t|x_1|^p - \gamma t(1-t) |x_0 - x_1|^p \quad \text{for every $t \in [0,1]$.}
\end{equation}

\begin{proof}[Proof of \cref{prop:Iconcavp}] Recall that we identify $\CC(\YY)$ with $\R^n$ and again for simplicity we prove the result only for $\beta =1$.

Therefore, letting $\psi_0$, $\psi_1 \in \R^n$, it holds
\begin{equation}
 \begin{split}
   - & (\psi_t)^{c,\epsilon}(x_t)  = \epsilon \log\sqa{ \sum\limits_{i=1}^n \exp\bra{\frac{(\psi_t)_i - |x_t - y_i|^p}{\epsilon}} \sigma_i  }  \\
 &\leq  \epsilon \log\Bigg[ \sum\limits_{i=1}^n \exp\bra{\frac{(1-t) (\psi_0)_i + t (\psi_1)_i -  (1-t) |x_0-y_i|^p - t|x_1-y_i|^p}{\epsilon}} \\
 & \quad   \cdot  \exp\bra{ \frac{ \gamma t(1-t) |x_0 - x_1|^p }{\epsilon}} \sigma_i \Bigg]
 \end{split}
\end{equation}
and we thus get, applying H\"older's inequality,
\begin{equation}
 (\psi_t)^{c,\epsilon}(x_t) \geq (1-t) (\psi_0)^{c,\epsilon}(x_0)  + t(\psi_1)^{c,\epsilon}(x_1) ) - t(1-t) \gamma  |x_0 - x_1|^p.
\end{equation}
 Therefore, taking the exponential and setting $h_t(x) = \exp((\psi_t)^{c,\eps}(x))$ we have
\begin{equation}
 h_t(x_t) \geq h_0(x_0)^{1-t} h_1(x_1)^t \exp(- t(1-t) \gamma |x_0 - x_1|^p).
\end{equation}

We next obtain from this bound an integral inequality, in a similar fashion to the Prékopa-Leindler inequality.  
Renormalizing, let us define $\hat{h}_k := h_k/ \|h_k\|_{L^1(\rho)}$ so that $\hat{h}_0$ and $\hat{h}_1$ are probability densities with respect to $\rho$. 

By \cref{lem:boundgeo} applied to the $p$-cost, $\rho_0 := \hat{h}_0 \rho$, $\rho_1 := \hat{h}_1 \rho$, using also the $\log$-concavity of $\rho$ and writing $\rho_t$ for the displacement interpolation between $\hat{h}_0\rho$ and $\hat{h}_1\rho$, $T$ the optimal map between $\hat{h}_0\rho$ and $\hat{h}_1\rho$, setting $T_t := (1-t) Id + tT$, we get for $\hat{h}_0\rho$-a.e.\ $x \in \XX$  the inequality :
\begin{equation}
 \begin{split}
  \rho_t(T_t(x)) &\leq (\hat{h}_0(x)\rho(x))^{1-t} (\hat{h}_1(T(x))\rho(T(x)))^t \\
&\leq \hat{h}_0(x)^{1-t} (\hat{h}_1(T(x)))^t \rho(x)^{1-t}\rho(T(x))^t \\
&\leq \hat{h}_0(x)^{1-t} (\hat{h}_1(T(x)))^t \rho(T_t(x))\\
&\leq h_t(T_t(x))\rho(T_t(x))   \exp(t(1-t) \gamma |x - T(x)|^p) \|h_0\|_{L^1(\rho)}^{-(1-t)} \|h_1\|_{L^1(\rho)}^{-t} .
 \end{split}
\end{equation}
so that 
\begin{equation}
\|h_0\|_{L^1(\rho)}^{1-t} \|h_1\|_{L^1(\rho)}^{t} \rho_t (x)  \exp(-t(1-t) \gamma |(T_t)^{-1}(x) - T((T_t)^{-1}(x))|^p) \leq  h_t(x)\rho(x)
\end{equation}
holds $\rho_t$ a.e., where the inverse of $T_t$ also makes sense $\rho_t$ a.e..
Integrating over  $\XX$ and using that $\rho_t = (T_t)_\sharp (h_0 \rho)$,
\begin{equation}
 \begin{split}
  \|h_t\|_{1} &\geq \|h_0\|_{L^1(\rho)}^{1-t} \|h_1\|_{L^1(\rho)}^{t} \int_{\XX} \exp(-\gamma t(1-t) |(T_t)^{-1}(x) - T((T_t)^{-1}(x))|^p)\rho_t(x)dx \\
&\geq \|h_0\|_{L^1(\rho)}^{1-t} \|h_1\|_{L^1(\rho)}^{t} \ds\int_{\XX} \exp(-\gamma t(1-t)|x - T(x) |^p) h_0(x) d\rho(x).
 \end{split}
\end{equation}
Applying Jensen's inequality with the exponential function we get,
\begin{equation}
 \|h_t\|_{L^1(\rho)} \geq \|h_0\|_{L^1(\rho)} ^{1- t} \|h_1\|_{L^1(\rho)} ^t \exp\bra{ - \gamma t(1-t) \ds\int_{\XX} |x - T(x) |^p h_0(x) d\rho(x)}.
\end{equation}
Finally, taking the $\log$ yields  the thesis:
\begin{equation}
 \begin{split}
  \log\I((1-t)\psi^0 + t \psi^1) &\geq  (1-t) \log\I(\psi^0) + t \log\I(\psi^1)\\
& \quad - t(1-t)\gamma \WW_p( \rho_{\psi_0}^{\eps}, \rho_{\psi_1}^{\eps})^p. \qedhere
 \end{split}
\end{equation}
\end{proof}
Using the established bound for $\log \I_\beta$, we prove the main result for this section.

\begin{proof}[Proof of \cref{theo:finalestp}]
From \cref{prop:Iconcavp},  writing $I_\beta := \log \I_\beta$ for brevity, one gets for all $t \in (0,1)$,
\begin{equation}
 \begin{split}
 \qquad \; \beta \gamma \WW_p( \rho_{\beta,\psi_0}^{\eps}, \rho_{\beta,\psi_1}^{\eps})^p &\geq \frac{(1-t) I_\beta(\psi^0) + t I_\beta(\psi^1) - I_\beta((1-t)\psi^0 + t \psi^1)}{t(1-t)} \\
&\geq - \int_0^1 \langle \nabla I_\beta(\psi_{st}) | \psi_1 - \psi_0 \rangle ds +\int_0^1 \langle \nabla I_\beta(\psi_{1 - (1-t)s}) | \psi_1 - \psi_0 \rangle ds.
 \end{split}
\end{equation}
Letting successively $t \to 0$ and $t \to 1$ we find
\begin{equation}
\beta \gamma  \WW_p( \rho_{\beta,\psi_0}^{\eps}, \rho_{\beta,\psi_1}^{\eps})^p \geq - \langle \nabla I_\beta(\psi_{0}) | \psi_1 - \psi_0 \rangle +\int_0^1 \langle \nabla I_\beta(\psi_{1 - s}) | \psi_1 - \psi_0 \rangle ds
\end{equation}
and 
\begin{equation}
 \beta\gamma  \WW_p( \rho_{\beta,\psi_0}^{\eps}, \rho_{\beta,\psi_1}^{\eps})^p \geq - \int_0^1 \langle \nabla I_\beta(\psi_{s}) | \psi_1 - \psi_0 \rangle ds + \langle \nabla I_\beta(\psi_{1}) | \psi_1 - \psi_0 \rangle
\end{equation}
so that, by summing the two inequalities,
\begin{equation}
2\beta\gamma  \WW_p( \rho_{\beta,\psi_0}^{\eps}, \rho_{\beta,\psi_1}^{\eps})^p \geq \langle \nabla I_\beta(\psi_1) - \nabla I_\beta(\psi_{0}) | \psi_1 - \psi_0 \rangle.
\end{equation}
Following the proof of \emph{Step 2} as proposed in \cref{prop:Kstrongconv} together with the additional comment in \cref{rem:sufficient_concav} and letting for $\phi_t^\eps = \psi_t^{c,\eps}$, for $t \in [0,1]$ where $\psi_t = (1-t) \psi_0 + t \psi_1$,
$M^\eps := \sup_{t \in [0,1]}\cur{\osc(\phi_t^\eps)}$,
and for all $\phi \in \CC(\XX)$, 
\begin{equation}\rho_{\phi} := \frac{e^{\phi}}{\int_\XX e^\phi d\rho} \rho, \end{equation}
we find
\begin{equation}\begin{split}\label{eq:finalestplambdadiscreteeps}
\beta e^{-\beta M^\eps }  \Var_\rho( \phi^{\eps}_1 - \phi^{\eps}_0) &  \leq e^{\beta M^\eps } \langle \nabla \KK^{c,\eps}(\psi_1) - \nabla \KK^{c,\eps}(\psi_0) | \psi_1 - \psi_0 \rangle \\
& \quad + 2\gamma \WW_p( \rho_{\beta \phi^{\eps}_0}, \rho_{\beta \phi^{\eps}_1})^p,
\end{split}\end{equation}
%
where we used \cref{eq:hessI}. Note that the scaling in $\beta^2$ of one term and in $\beta$ of the other in \cref{eq:hessI} explains the presence of $\beta$ on the right hand side of \cref{eq:finalestplambdadiscreteeps}.

If we then let $\eps \downarrow 0$ minding \cref{rem:vanishing-entropy}, by choosing for simplicity the potentials so that $\phi_k$ have zero $\rho$-mean, we get that for all $\beta >0$ the following bound holds:
\begin{equation}\label{eq:finalestplambda}
\beta e^{-\beta M_{\phi_0,\phi_1}}  \| \phi_1 - \phi_0 \|_{L^2(\rho)}^2  \leq e^{\beta M_{\phi_0,\phi_1}} \langle \mu_0 - \mu_1 | \psi_1 - \psi_0 \rangle + 2\gamma \WW_p( \rho_{\beta \phi_0}, \rho_{\beta \phi_1})^p.
\end{equation}

The next step is to derive an estimate of $\WW_p(\rho_{\beta \phi_0}, \rho_{\beta \phi_1})$ in terms of $\nor{\phi_1-\phi_0}_{L^2(\rho)}$. Again, for simplicity we argue in the case $\beta = 1$. First, we notice that by the monotonicity of Wasserstein distances, we have $\WW_p(\rho_{\phi_0}, \rho_{\phi_1}) \le \WW_2(\rho_{\phi_0}, \rho_{\phi_1})$.
Then,  we apply \cref{lem:peyre}, with $\mu_k = \rho_{\phi_k}$, $k \in \cur{0,1}$, so that $f_k \geq e^{-M_{\phi_0,\phi_1}}$ and $\rho$ is $\log$-concave with bounded support so it has finite Poincaré constant. We get

\begin{equation}\label{eq:W2boundL2}
\WW_2^2(\rho_{\phi_1}, \rho_{\phi_0}) \leq C_\rho e^{M_{\phi_0,\phi_1}} \nor{  \frac{e^{\phi_1}}{\int_\XX e^{\phi_1} d \rho}  - \frac{e^{\phi_0}}{\int_\XX e^{\phi_0} d\rho}  }_{L^2(\rho)}^2.
\end{equation}
For any $x \in \XX$, we have 
\begin{equation}
 \begin{split}
  \abs{\frac{e^{\phi_1(x)}}{\int_\XX e^{\phi_1} d\rho} - \frac{e^{\phi_0(x)}}{\int_\XX e^{\phi_0} d\rho} } &\leq \abs{\frac{e^{\phi_1(x)} - e^{\phi_0(x)} }{\int_\XX e^{\phi_1} d\rho}} + e^{\phi_0(x)}\abs{\frac{1}{\int_\XX e^{\phi_1} d\rho} - \frac{1}{\int_\XX e^{\phi_0} d\rho} } \\
&\leq e^{M_{\phi_0,\phi_1}} | e^{\phi_1(x)} - e^{\phi_0(x)}| + e^{2M_{\phi_0,\phi_1}} \abs{ \int_\XX  \bra{e^{\phi_1} - e^{\phi_0}} d\rho }\\
& \leq  e^{2M_{\phi_0,\phi_1}} \Big( |\phi_1(x) - \phi_0(x)| +  \|\phi_1 - \phi_0\|_{L^1(\rho)} \Big),
 \end{split}
\end{equation}
where the zero $\rho$-mean property of $\phi_k$ assures that
\begin{equation} -M_{\phi_0,\phi_1} \leq \min\phi_k \leq \max\phi_k\leq M_{\phi_0,\phi_1}.\end{equation}
Integrating with respect to $\rho$ yields
\begin{equation}
 \nor{  \frac{e^{\phi_1}}{\int_\XX e^{\phi_1} d \rho}  - \frac{e^{\phi_0}}{\int_\XX e^{\phi_0} d\rho}  }_{L^2(\rho)}^2 \le C_\XX e^{4M_{\phi_0,\phi_1} } \nor{ \phi_1 - \phi_0}_{L^2(\rho)}^2.
\end{equation}
From this estimation, \cref{eq:finalestplambda} and \cref{eq:W2boundL2} we get
\begin{equation}
\beta e^{-\beta M_{\phi_0,\phi_1}}  \| \phi_1 - \phi_0 \|_{L^2(\rho)}^2  \leq e^{\beta M_{\phi_0,\phi_1}} \langle \mu_0 - \mu_1 | \psi_1 - \psi_0 \rangle + \gamma C_\rho e^{\frac{5p\beta}{2}M_{\phi_0,\phi_1}} \beta^p \|\phi_1 - \phi_0\|_{L^2(\rho)}^p
\end{equation}
so,
\begin{equation}
\beta e^{-2\beta M_{\phi_0,\phi_1}}  \| \phi_1 - \phi_0 \|_{L^2(\rho)}^2 -  \gamma C_\rho e^{\frac{(5p-2)\beta}{2}M_{\phi_0,\phi_1}} \beta^p \|\phi_1 - \phi_0\|_{L^2(\rho)}^p  \leq \langle \mu_0 - \mu_1 | \psi_1 - \psi_0 \rangle.
\end{equation}

If we forget for a moment the exponential terms, the optimal choice is 
\begin{equation}
\beta = \frac{ \|\phi_1 - \phi_0\|_{L^2(\rho)}^\frac{2-p}{p-1}}{(\gamma C_\rho)^\frac{1}{p-1}},
\end{equation} 
and this is the scaling that guides us in the following calculations. We now set
$C =  \gamma C_\rho$, $\alpha = 2M_{\phi_0,\phi_1}$ and $\alpha' = \frac{(5p-2)}{2}M_{\phi_0,\phi_1} $ so that $C, \alpha, \alpha' > 0$.

Letting
\begin{equation}
\begin{split}
f(\beta) &:= \beta e^{-\alpha \beta}  || \phi_1 - \phi_0 ||_{L^2(\rho)}^2 -  C e^{\alpha'\beta} \beta^p ||\phi_1 - \phi_0||_{L^2(\rho)}^p \\
g(\eps) &:= \eps e^{-\alpha \eps  ||\phi_1 - \phi_0||_{L^2(\rho)}^\frac{2-p}{p-1} } -  C e^{\alpha' \eps ||\phi_1 - \phi_0||_{L^2(\rho)}^\frac{2-p}{p-1} } \eps^p \\
h(\beta) &:= \beta e^{-\alpha \beta}  -  C e^{\alpha'\beta} \beta^p 
\end{split}
\end{equation}

so that 
\begin{equation}
f(\beta) = ||\phi_1 - \phi_0||_{L^2(\rho)}^\frac{p}{p-1} g(\eps) \, \text{ with $\beta = \eps  ||\phi_1 - \phi_0||_{L^2(\rho)}^\frac{2-p}{p-1}$ }.
\end{equation}
Noticing
\begin{equation}
\left\{
    \begin{array}{ll}
        g(\eps) \geq h(\eps)& \mbox{if  } ||\phi_1 - \phi_0||_{L^2(\rho)} \leq 1 \\
        f(\beta) \geq h(\beta) ||\phi_1 - \phi_0||_{L^2(\rho)}^2 & \mbox{if } ||\phi_1 - \phi_0||_{L^2(\rho)} \geq 1
    \end{array}
\right.
\end{equation}
we get


\begin{equation}\label{eq:uglyccl}
\min (\| \phi_1 - \phi_0 \|_{L^2(\rho)}^2, \|\phi_1 - \phi_0\|_{L^2(\rho)}^\frac{p}{p-1})  \leq \frac{1}{\sup\limits_{\beta \in \R_+} h(\beta) } \langle \mu_0 - \mu_1 | \psi_1 - \psi_0 \rangle.
\end{equation}
Now, deriving the absolute bound
\begin{equation}
\begin{split}
\| \phi_1 - \phi_0 \|_{L^2(\rho)} &\le |\XX| (\| \phi_0 \|_{L^\infty} + \| \phi_1 \|_{L^\infty}) \\
& \le 2 |\XX| M_{\phi_0, \phi_1} \le C = C(\rho, p, \YY)
\end{split}
\end{equation}
we get $ \|\phi_1 - \phi_0\|_{L^2(\rho)}^\frac{p}{p-1} \le C \| \phi_1 - \phi_0 \|_{L^2(\rho)}^2$.
From this last inequality and \cref{eq:uglyccl} we can conclude
\begin{equation}
 \|\phi_1 - \phi_0\|_{L^2(\rho)}^\frac{p}{p-1} \le C \langle \mu_0 - \mu_1 | \psi_1 - \psi_0 \rangle
\end{equation}
for some $C = C(p, \rho, \XX, \YY) > 0$.
One final application of \cref{lem:approx} extends the bounds from the semi-discrete case to the general case of compact $\YY$. \qedhere

\end{proof}

\begin{rem} One can easily check that the entire proof above relies only on (Lipschitz continuity and) the curvature condition \cref{eq:pPC} satisfied by the $p$-cost. Thus, if the cost $c$ is such that, for some $p>1$ and $\gamma<\infty$, it holds
\begin{equation}\label{eq:cpPC}
c(x_t,y) \geq (1-t) c(x_0,y) + tc(x_1,y) - \gamma t(1-t) \abs{x_0 - x_1}^p \quad \text{ for all $x_0, x_1 \in \XX$ and $y \in \YY$}
\end{equation}
then \cref{prop:Iconcavp} still holds. We notice that such property, which we may refer to as a  $(p,\gamma)$-curvature condition for the cost function $c$, is stable by taking sums and infima. As an example, it holds for the boundary cost analogue to \cref{eq:boundary-2} with general exponent $p$.
\end{rem}

\begin{rem}\label{rem:density-p}
Again, we notice that the  only point where $\log$-concavity of $\rho$ is used is in the proof of \cref{prop:Iconcavp}. Therefore, in order to cover the case of $\rho$ that admits a density uniformly bounded from above and below by strictly positive constants with respect to a $\log$-concave measure $\lambda$, we need to replace $\rho$ with $\lambda$ in the definition of $\I$. This leads to a version of \cref{eq:logconcavIp} with $\lambda^\eps_{\beta, \psi_k}$ instead of $\rho^{\eps}_{\beta, \psi_k}$. However, this changes very little in the proof of \cref{theo:finalestp}: we obtain \cref{eq:finalestplambdadiscreteeps} arguing as in \cref{prop:Kstrongconv} and taking into account \cref{rem:density}; subsequently, we apply \cref{lem:peyre} with respect to the measure $\lambda$; and finally we replace $L^2(\lambda)$ with $L^2(\rho)$, up to a multiplicative constant.
 \end{rem}

\section{Quantitative stability of maps}\label{sec:maps}

In this section, we focus on the stability of optimal transport maps by establishing \cref{thm:main-map}. We make the following standing assumption: $\rho$ is a probability measure with bounded and convex support $\XX \subseteq\R^d$, with density bounded below and above by two constants:
\begin{equation}\label{eq:density-rho}
 0<m_\rho \le \rho(x) \le  M_\rho < +\infty \quad \text{for every $x \in \XX$.}
\end{equation}

Our goal is now to use the results from the previous section to yield a bound on the $L^2$ distance between two optimal maps $T_{0}$ and $T_{1}$ from the source $\rho$ to two target distributions $\mu_0$, $\mu_1$, supported on a bounded set $\YY$. By the extension of Brenier theorem to general $p$-power costs \cite[Theorem 6.2.4]{AGS} ($p>1$), the optimal transport map $T$ for $\WW_p(\rho, \mu)$ is given by
\begin{equation}\label{eq:brenier-gradient}
 T(x) = x - (\nabla \phi(x))^{(q/p)} \quad \text{for a.e.\ $x \in \XX$,}
\end{equation}
where $\phi$ is the Kantorovich potential, $q=p/(p-1)$ and we recall the notation $v^{(\alpha)} = |v|^{\alpha-1} v$. Thus, our strategy is to establish first upper bounds on $\|\nabla \phi_1 - \nabla \phi_0 \|_{L^2}$ in terms of $\|\phi_1- \phi_0\|_{L^2}$,  whenever  $\phi_k$, $k \in \cur{0,1}$ are $c$-concave functions. In the following section, we derive such bounds, first for the quadratic case -- recalling the argument from \cite{delalande2022quantitative} -- and then more generally for $c$-concave functions. 

\subsection{The quadratic case}

We report in this section the argument from \cite{delalande2022quantitative}. Let us recall that the Brenier map is the optimal transport map for the $p$-cost when $p=2$.

\begin{theo}
Let $\XX$, $\YY$ be compact subsets of $\R^d$ with $\XX$ convex, let $\rho$ be a probability measure with density w.r.t.\ the Lebesgue measure bounded below and above by strictly positive constants as in \cref{eq:density-rho} and  let $\mu_0$, $\mu_1 \in \PP(\YY)$. Then, there exists a constant $C = C(d, \XX,\YY, m_\rho, M_\rho)< \infty$ such that
\begin{equation}
\|T_0 - T_1 \|_{L^2(\rho)} \leq C \WW_1(\mu_0, \mu_1)^\frac{1}{6}.
\end{equation}
where $T_k$ denotes the Brenier map from $\rho$ to $\mu_k$, for $k \in \cur{0,1}$.
\end{theo}

This theorem is an immediate consequence of \cref{theo:stabcompact}, \cref{eq:brenier-gradient} which reads
\begin{equation}
 T(x) =  x - \nabla \phi(x) = \nabla \bra{ \frac {|x|^2}{2} - \phi(x)}, \quad \text{for a.e.\ $x \in \XX$,}
\end{equation}
so that $T(x)$ is the gradient of the function $x \mapsto |x|^2/2-\phi(x)$ which is convex and Lipschitz, and the following ``reverse'' Poincaré inequality,  \cite[Proposition 5.11]{delalande2022quantitative} -- but see also \cite{shashiashvili2005estimation, hussain2008weighted, pevcaric2015reverse} for similar results -- whose statement we report here for the sake of completeness and for comparison with our extension to fractional spaces, \cref{prop:p-convgaliardo} below.

\begin{prop}\label{prop:convgaliardo} Let $\XX\subseteq \R^d$ be a bounded domain with rectifiable boundary and let $u,v: \XX \to \R$ be convex Lipschitz functions. Then, there exists a constant $C = C(d)<\infty$ depending on $d$ only such that
\begin{equation}
\| \nabla u - \nabla v\|_{L^2(\XX)}^2 \leq C \HH^{d-1}(\partial \XX)^\frac{2}{3} ( \|\nabla u \|_{L^\infty(\XX)} + \|\nabla v \|_{L^\infty(\XX)})^\frac{4}{3} \|u-v\|_{L^2(\XX)}^\frac{2}{3},
\end{equation}
where $\HH^{d-1}$ denotes the $(d-1)$-dimensional Hausdorff measure. 
\end{prop}

\subsection{The case $p>2$}

This is a minor modification of the quadratic case, keeping in mind that the optimal transport map depends on the gradient of the Kantorovich potential in a H\"older continuous way. We prove the following result.

\begin{theo}\label{theo:stabmapgeq2} Let $\XX$, $\YY$ be compact subsets of $\R^d$ with $\XX$ convex, let $\rho$ be a probability measure, with density w.r.t.\ the Lebesgue measure bounded below and above by strictly positive constants as in \cref{eq:density-rho}. For $k \in \cur{0,1}$, let $\mu_k \in \PP(\YY)$ and denote with $T_{k}$ the optimal transport map for $\WW_p(\rho, \mu_k)$. Then, there exists a constant $C = C(p, d, \XX, \YY, m_\rho, M_\rho)< \infty$  such that
\begin{equation}
\|T_0 - T_1 \|_{L^2(\rho)} \leq C \WW_1(\mu_0, \mu_1)^\frac{1}{6(p-1)}.
\end{equation}
\end{theo}

\begin{proof}
We start from the inequality $|v^{(\alpha)} - w^{(\alpha)}| \leq C(\alpha) |v-w|^{\alpha}$, valid for all $0<\alpha<1$ and $v,w \in \R^d$
with $C = C(\alpha,d)<\infty$. We find
\begin{equation}
\begin{split}
 \nor{ T_0 - T_1}_{L^2(\rho)}^2 & = \int_\XX |\nabla \phi_0^{(q/p)} - \nabla \phi_1^{(q/p)} |^2 d\rho \\
& \leq  C \int_\XX | \nabla \phi_0 - \nabla \phi_1|^{2q/p} d\rho \\
&\leq C \nor{  \nabla \phi_0 - \nabla \phi_1}_{L^2(\rho)}^{2q/p},
\end{split}
\end{equation} 
having used H\"older inequality and the fact that $q/p= 1/(p-1)<1$ for $p>2$. Next, we notice that there exists $\gamma=\gamma(p, \XX, \YY)<\infty$ such that, for any function $\phi$ that is $c$-concave with respect to the $p$-cost, the function  $x \mapsto \phi (x)- \gamma |x|^2$ is concave and Lipschitz (with a  bound on the Lipschitz constant uniform on $\phi$). We then proceed with the upper bounds (where the constant $C$ may change from line to line)
\begin{equation}
\begin{split}
 \nor{  \nabla \phi_0 - \nabla \phi_1}_{L^2(\rho)}^{2} & \le C \nor{  \nabla \phi_0 - \nabla \phi_1}_{L^2(\XX)}^{2}\\
 & = C  \nor{  \nabla (\phi_0 - \gamma|\cdot|^2) - \nabla(\phi_1 - \gamma|\cdot|^2)}_{L^2(\XX)}^{2}\\
 & \le C  \nor{ \phi_0-\phi_1}_{L^2(\XX)}^\frac{2}{3}  \text{by \cref{prop:convgaliardo}}\\
 & \le C \nor{ \phi_0-\phi_1}_{L^2(\rho)}^\frac{2}{3} \\
 & \le C \WW_1(\mu_0, \mu_1)^\frac{1}{3}.
\end{split}
\end{equation}
where the last inequality follows by \cref{thm:main-pot} assuming without loss of generality that $\phi_0$, $\phi_1$ have zero $\rho$-mean.
\end{proof}

\subsection{The case $1<p<2$}

This case is instead less obvious, again because there is no shift reducing a $c$-concave function to a concave function when $c$ is the $p$-cost for $1<p<2$. We prove the following result.

\begin{theo}\label{theo:stabmapleq2}
Let $p \in (1,2)$, set $q=p/(p-1)$, let $\XX$, $\YY \subseteq \R^d$ be compact subsets with $\XX$ convex, let $\rho$ be a $\log$-concave probability measure on $\XX$ such that \cref{eq:density-rho} holds, let $\mu_0$, $\mu_1 \in \PP(\YY)$ and denote with $T_{k}$ the optimal transport map for $\WW_p(\rho, \mu_k)$.  Then, for every $\theta \in (0,(p-1)/(p+1))$ there exists a constant $C = C(\theta, p, d, \XX, \YY, m_\rho, M_\rho)< \infty$  such that
\begin{equation}
\|T_0 - T_1 \|_{L^2(\rho)} \leq C  \WW_1(\mu_0,\mu_1)^\frac{\theta}{q}. 
\end{equation}
\end{theo}

To establish this theorem we need an equivalent of \cref{prop:convgaliardo} for $c$-concave functions with respect to the $p$-cost.

\begin{prop}\label{prop:p-convgaliardo} Let $p \in (1,2)$, $\XX \subseteq \R^d$ be a bounded domain with rectifiable boundary and let $u,v: \XX \to \R$ be $c$-concave with respect to the $p$-cost. Then, for all $\theta \in (0, (p-1)/(p+1))$ there exists a constant $C = C(\theta,\XX, p)<\infty$ such that
\begin{equation}
\| \nabla u - \nabla v\|_{L^2(\XX)} \leq C \|u-v\|_{L^2(\XX)}^{\theta}.
\end{equation}
\end{prop}

For comparison with our argument, let us recall that the proof of \cref{prop:convgaliardo} relies on two ingredients. First, one establishes the one-dimensional inequality \cite[Proposition 5.24]{delalande2022quantitative}
\begin{equation}\label{eq:convgagliardo}
\int_0^1 |u' - v'|^2 \leq 8\bra{ \operatorname{Lip}(u) + \operatorname{Lip}(v) }^{\frac{4}{3}}\bra{  \int_0^1 |u - v|^2}^{\frac 1 3}
\end{equation}
that settles the case $\XX =[0,1]\subseteq \R$ (and by scaling any other interval). Then, for the general case, one applies an integral-geometric argument \cite[Lemma 5.25]{delalande2022quantitative} decomposing the $d$-dimensional integrals  into integral over segments. In our context, the integral-geometric argument formula will be identical, up to harmless variations in the exponents. Our aim is therefore to establish a suitable analogue of \cref{eq:convgagliardo} for functions that are $c$-concave with respect to the $p$-cost. Before we do this, it is convenient to notice the following consequence of $c$-concavity with respect to the $p$-cost. 

\begin{lem}\label{lem:NCpconc}
Let $p \in (1,2)$, $\XX \subseteq \R^d$ be open and bounded and $\phi:\XX \to \R$ be $c$-concave with respect to the $p$-cost. Then, at every  $x$, $y \in \XX$ that are differentiability points for $\phi$, it holds
\begin{equation}\label{eq:plambdaconc}
\langle \nabla \phi(x) - \nabla \phi(y),  x-y \rangle \leq \frac{2\gamma}{p^2} |x-y|^p,
\end{equation}
where $\gamma = \gamma(p)$ is the same constant as in \cref{eq:pPC}.
\end{lem}

\begin{proof}  By $c$-concavity, the graph of the associated transport map
\begin{equation}
 x \to T(x) := x - (\nabla \phi(x))^{(q/p)}
\end{equation}
where $q=p/(p-1)$, that is defined at every differentiability point of $x$, is $c$-cyclically monotone, i.e.,
\begin{equation}\label{eq:c-monotonicity-T}
 |T(x) - x|^p - |T(x) - y|^p \leq  |T(y) - x|^p-|T(y) - y|^p.
\end{equation}
Next, we recall the inequality, valid for any $v,h,z \in \R^d$,
\begin{equation}\label{eq:taylor-inequality}
 \begin{split}
  |v+h-z|^p &- |v-z|^p - p\langle (v-z)^{(p-1)}, h \rangle \\
& = p\int_0^1 \langle (v +th -z)^{(p-1)} - (v - z)^{(p-1)}, h \rangle dt \\
& \leq \gamma \int_0^1 t^{p-1} |h|^p dt = \frac{\gamma}{p} |h|^p,
 \end{split}
\end{equation}
where we used the inequality \cref{eq:tooltocurvature}. 
Applying \cref{eq:taylor-inequality} in  both sides of \cref{eq:c-monotonicity-T}, i.e., first with $z = T(x)$, $h = y-x$ and then for $z = T(y)$, $h = x-y$, we get
\begin{equation}
 -p\langle (x-T(x))^{(p-1)}, y-x \rangle - \frac{\gamma}{p} |y-x|^p \leq p \langle (y-T(y))^{(p-1)}, x-y \rangle + \frac{\gamma}{p} |y-x|^p
\end{equation}
so that, noticing that $(x-T(x))^{(p-1)} = \nabla \phi(x)$ by the very definition of $T$, we obtain the thesis.
\end{proof}

 In the proof below we only use that $c$-concave functions for the $p$-cost are Lipschitz (with uniform bounds on their Lipschitz constant) and satisfy inequality \cref{eq:plambdaconc}. Thus, we introduce the following notation: for $p \in (1,2)$,  $\lambda \in \R$, we say that a Lipschitz function $u: \XX \to \R$ is $(p,\lambda)$-concave if, for every $x$, $y$ that are differentiability point of $u$ it holds
\begin{equation}\label{eq:weak-hess-neg}
\langle \nabla u(x) - \nabla u(y) ,  x-y \rangle \leq \lambda |x-y|^p.
 \end{equation}
We are now in a position to state and prove the following analogue of \cref{eq:convgagliardo}.

\begin{lem}\label{lem:pconcavgagliardo} Let $p \in (1,2)$ and $\lambda, L\ge 0$. For every $\theta \in (0,(p-1)/(p+1))$, there exists $C= C(p,\lambda, L, \theta)<\infty$ such that for $u,v: [0,1] \to \R$ both $(p,\lambda)$-concave and $L$-Lipschitz,  it holds
\begin{equation}\label{eq:p-convgagliardo}
\int_0^1 |u' - v'|^2 dx \leq C\bra{  \int_0^1 |u - v|^2dx}^{\theta}.
\end{equation}
\end{lem}

\begin{proof}[Proof of \cref{lem:pconcavgagliardo}] 

Without loss of generality, we may assume throughout the proof that $\int_0^1 u dx = \int_0^1 v dx = 0$.

For comparison with the proof in our case, we recall that a key point in the proof of \cref{eq:convgagliardo} in \cite{delalande2022quantitative} is the following  integration by parts (assuming by approximation that the functions are $\CC^2$ smooth):
\begin{equation}
 \begin{split}
 \int_0^1 |u' - v'|^2  &= (u-v)(u'-v') \Big|_0^1 - \int_0^1 (u-v)(u''-v'') \\
&\leq 2(\Lip(u) + \Lip(v) ) \|u-v\|_{\infty} + \|u-v\|_{\infty} \int_0^1 (|u''| + |v''|) \\
& = 2(\Lip(u) + \Lip(v) ) \|u-v\|_{\infty} + \|u-v\|_{\infty} \int_0^1 (u'' + v'')\\
&\leq \sqa{ 2(\|u'\|_{L^\infty} + \|v'\|_{L^\infty})  +  (u' + v')\Big|_0^1} \|u-v\|_{\infty} \\
&\leq 4(\|u'\|_{L^\infty} + \|v'\|_{L^\infty}) \|u-v\|_{\infty},
 \end{split}
\end{equation}
where the convexity assumption on $u$ and $v$ allowed us to remove the absolute values in the terms $|u''|$ and $|v''|$.
  The above argument can be interpreted as an interpolation bound on the $H^{1,2}$ norm in terms of the $L^\infty$ norm of the function and the $L^1$ norm of its second derivative, that together with the positivity condition reduces  to a $L^\infty$ norm of the derivative. Finally, the standard Gagliardo-Nirenberg inequality allows to replace the $L^\infty$ norm in the right hand side with the $L^2$ norm (which gives the exponents $4/3$ and $2/3$).

  In our case of $(p,\lambda)$-concave functions, the positivity of second derivatives is replaced by the weaker condition \cref{eq:weak-hess-neg}, so the idea is to interpolate with fractional derivatives of order between $1$ and $2$. We start from the assumption, which in the one dimensional case reads
  \begin{equation}
   (u'(x)-u'(y))(x-y) - \lambda |x-y|^p \leq 0
  \end{equation}
  for every  $x, y \in (0,1)$ where  $u$ is differentiable, which is a set of full Lebesgue measure, because $u$ is Lipschitz continuous. Without loss of generality, we may also assume that $\lambda \ge 0$. Therefore, we have the inequality (assuming $x \neq y$):
\begin{equation}\begin{split}
         \frac{|u'(x) - u'(y)|}{|x-y|} &= \frac{|(u'(x) - u'(y))(x-y)|}{|x-y|^2} \\
&\leq \frac{|(u'(x) - u'(y))(x-y) -  \lambda |x-y|^p| + \lambda |x-y|^p}{|x-y|^2}\\
&= \frac{(u'(y) - u'(x))}{x-y} +  \frac{2 \lambda}{|x-y|^{2-p}}.
                \end{split}
\end{equation}
Given $\alpha \in (0,1)$, we bound from above
\begin{equation}\begin{split}
[u']_{\alpha,1} & = \int_0^1 \int_0^1 \frac{|u'(x) - u'(y)|}{|x-y|^{1+\alpha}} dydx\\
& \leq \int_0^1 \int_0^1 \frac{u'(y) - u'(x)}{(x-y) |x-y|^\alpha} dydx +  \int_0^1 \int_0^1 \frac{2\lambda dxdy}{|x-y|^{2-p+\alpha}} \\
&= - \int_0^1 \int_{-1}^1 \indicatrice_{\cur{0 \leq x +h \leq 1} }\frac{u'(x+h) - u'(x)}{h |h|^\alpha} dhdx + 2\lambda \int_0^1 \int_{-1}^1 \indicatrice_{\cur{0 \leq x +h \leq 1}} \frac{dhdx}{|h|^{2-p+\alpha}}.
\end{split}
\end{equation}
where $\chi$ denotes the indicator function. We bound separately the two terms. For the second one, we simply use the fact that $\int_0^1 h^{-\beta} ds<\infty$ if $\beta<1$, hence
\begin{equation}
 2\lambda \ds\int_0^1 \int_{-1}^1 \indicatrice_{\cur{0 \leq x +h \leq 1}} \frac{dhdx}{|h|^{2-p+\alpha}} \le C,
\end{equation}
for some $C = C(p,\alpha, \lambda)<\infty$, provided that $2-p+\alpha<1$, i.e., $\alpha<p-1$. For the first term, we integrate first with respect to $x \in [0,1]$, so that for fixed $h \in [-1,1]$, we find
\begin{equation}\begin{split}
\int_0^1 \indicatrice_{\cur{0 \leq x +h \leq 1} } \bra{u'(x+h) - u'(x)} dx&  = u(\min\cur{1+h, 1 }) - u(\max\cur{h,0 } ) \\
& \quad  -  u(\min\cur{1-h, 1 }) + u(\max\cur{-h,0}).
\end{split}
 \end{equation}
 When we integrate also with respect to $h \in [-1,1]$ it is convenient to consider separately the two cases $h>0$ and $h<0$. In the first case, we find
 \begin{equation}
 \begin{split}
    \int_{0}^1 \int_0^1 \indicatrice_{\cur{0 \leq x +h \leq 1} }\frac{u'(x+h) - u'(x)}{|h|^\alpha h} dxdh & =  \int_{0}^1 \frac{ u(1) -  u(1-h)- u(h)  + u(0)}{h^{1+\alpha}} dh \\
   & \le 2 \Lip (u) \int_0^1 \frac{1}{h^\alpha} dh  \le C,
\end{split}
 \end{equation}
 where $C = C(\alpha, L)<\infty$, having used that $\alpha <1$. In the case $h<0$, we find similarly
 \begin{equation}
 \begin{split}
  &  \int_{-1}^0 \int_0^1 \indicatrice_{\cur{0 \leq x +h \leq 1} }\frac{u'(x+h) - u'(x)}{|h|^\alpha h} dxdh \\
  = \quad & \int_{-1}^0 \frac{ u(1+h) -  u(0)- u(1)  + u(-h)}{|h|^\alpha h} dh \\
  \le \quad  & 2 \Lip (u) \int_0^1 \frac{1}{h^\alpha} dh  = C.
\end{split}
 \end{equation}
 Thus, we collect the bound
 \begin{equation}
  [u']_{\alpha,1} \le C = C(p, \alpha, \lambda, L)<\infty,
 \end{equation}
 provided that $0<\alpha<p-1$. Since $u$ is Lipschitz (and has zero mean), it also holds
 \begin{equation}
   \nor{u}_{W^{1,1}((0,1))} \le  \nor{u}_{W^{1,\infty}((0,1))} \le C.
 \end{equation}
In conclusion, the thesis follows as an application of \cref{eq:gagliardo-useful} on $\Omega = (0,1)$ and $r := 1+\alpha \in (1, p)$ and letting $\theta :=1-2/(1+r)$, which can be any real in the interval $(0, (p-1)/(p+1))$.
\end{proof}

We now indicate how to complete the proof of \cref{prop:p-convgaliardo} following the steps of  \cite[Proposition 5.11]{delalande2022quantitative}, for which we refer for more details.

\begin{proof}[Proof of \cref{prop:p-convgaliardo}.]

The first step is to notice the good scaling properties of inequality \cref{eq:p-convgagliardo}. Through an affine change of variable, one easily proves that if $p \in (1,2)$, $\lambda, L, l \ge 0$ and $\theta \in (0,(p-1)/(p+1))$, there exists $C= C(p,\lambda, L, l, \theta)<\infty$ such that for any interval $I \subset \R$ with $|I| \le l$, $u,v: I \to \R$ both $(p,\lambda)$-concave and $L$-Lipschitz,  it holds
\begin{equation}\label{eq:p-convgagliardo-scaled}
\int_I |u' - v'|^2 dx \leq C\bra{  \int_I |u - v|^2 dx}^{\theta}.
\end{equation}
To pursue the proof we need to introduce some new notation that we try to keep as light as possible, referring to \cite{delalande2022quantitative} and the references therein for more details.

Let $\LL^d$ be the set of affine lines in $\R^d$ identified with $\{ \ell = (e,p) \in \R^d \times \R^d | e \in \S^{d-1}, p \in e^\perp\}$ equipped with the Riemannian metric. We write $d\LL^d$ the volume measure, so that for any $\phi : \LL^d \to \R$
\begin{equation}
\ds\int_{\ell \in \LL^d} \phi(\ell)d\LL^d(\ell) = \int_{e \in \S^{d-1}} \int_{p \in e^\perp} \phi((e,p))d\mathcal{H}^{d-1}(p)d\mathcal{H}^{d-1}(e)
\end{equation}
where $\mathcal{H}^{d-1}$ denotes the $(d-1)$-dimensional Hausdorff measure. 

Back to our setting, by \cite[Lemma 5.25]{delalande2022quantitative}, one obtains
\begin{equation}\label{eq:int-form-u}
\| u \|_{L^2(\XX)}^2 = C \ds\int_{\ell \in \LL^d} \int_{y \in \ell \cap \XX} u(y)^2 dy d\LL^d(\ell)
\end{equation}
and 
\begin{equation}\label{eq:int-form-grad-u}
\| \nabla u \|_{L^2(\XX)}^2 = C' \ds\int_{\ell \in \LL^d} \int_{y \in \ell \cap \XX} \la \nabla u(y), e(\ell) \ra ^2 dy d\LL^d(\ell)
\end{equation}
where $C = C(d),  C'=C'(d) \in (0,\infty)$ are dimensional constants. Since $\XX$ is convex, we can write $\ell \cap \XX =  I_\ell$  for a suitable line segment  (that we identify with a real interval).
Then we use \cref{eq:p-convgagliardo-scaled} with $t \in I_\ell \to u(p + te)$ and $l = \diam(\XX)$ for a.e. $\ell \in \LL^d$, \cref{eq:int-form-u}, \cref{eq:int-form-grad-u} and H\"older's inequality obtaining the thesis, \begin{equation}
\| \nabla u - \nabla v\|_{L^2(\XX)} \leq C \|u-v\|_{L^2(\XX)}^{\theta}.\qedhere
\end{equation}
\end{proof}

We end this section with the proof of \cref{theo:stabmapleq2}, which in turn concludes the proof of \cref{thm:main-map}.

\begin{proof}[Proof of \cref{theo:stabmapleq2}]
Naming $\phi_0$ and $\phi_1$ the zero-mean (with respect to $\rho$) Kantorovich potentials associated to $T_0$ and $T_1$, by \cref{eq:brenier-gradient} we have
\begin{equation}
 \| T_0 - T_1 \|_{L^2(\XX)} = \nor{ (\nabla \phi_0)^{(q/p)} - (\nabla \phi_1)^{(q/p)} }_{L^2(\XX)},
\end{equation}
where $q= p/(p-1)$. Using the inequality
\begin{equation}
 |v^{(\alpha)} - w^{(\alpha)}| \leq C |v-w|( |v|^{\alpha-1}+ |w|^{\alpha-1}),
\end{equation}
valid for $v$, $w \in \R^d$, $\alpha >1$ with $C = C(\alpha,d)<\infty$, we find for $\alpha:=q/p>1$,
\begin{equation}
  \| T_0 - T_1\|_{L^2(\XX)}  \le C \bra{ \Lip(\phi_0)^{q/p-1} + \Lip(\phi_1)^{q/p-1} } \nor{ \nabla \phi_0 - \nabla \phi_1}_{L^2(\XX)}.
\end{equation}
Since the Lipschitz constant of $\phi_0$ and $\phi_1$ are uniformly bounded by $c$-concavity,  it follows that
\begin{equation}\begin{split}
 \| T_0 - T_1\|_{L^2(\XX)}  & \le C  \nor{ \nabla \phi_0 - \nabla \phi_1}_{L^2(\XX)}\\
 &\le C \nor{ \phi_0 - \phi_1}_{L^2(\XX)}^{\theta}\\
 & \le C  \nor{ \phi_0 - \phi_1}_{L^2(\rho)}^{\theta}
 \end{split}
\end{equation}
where the constant $C$ may change from line to line, and we applied \cref{prop:p-convgaliardo}. To conclude, we apply \cref{theo:finalestp} and obtain
\begin{equation}
  \| T_0 - T_1\|_{L^2(\rho)} \le C \WW_1(\mu_0, \mu_1)^{\frac{ \theta }{q } },
\end{equation}
which is the thesis.
\end{proof}

\section*{Acknowledgements and Funding} Both authors thank A.\ Pratelli, L.\ Ambrosio, S.\ Di Marino and A.\ Gerolin for useful discussions and comments during the development of this work.

D.T.\ acknowledges the MUR Excellence Department Project awarded to the Department of Mathematics, University of Pisa, CUP I57G22000700001,  the HPC Italian National Centre for HPC, Big Data and Quantum Computing - Proposal code CN1 CN00000013, CUP I53C22000690001, the PRIN 2022 Italian grant 2022WHZ5XH - ``understanding the LEarning process of QUantum Neural networks (LeQun)'', CUP J53D23003890006, the INdAM-GNAMPA project 2024 ``Tecniche analitiche e probabilistiche in informazione quantistica'' and the project  G24-202 ``Variational methods for geometric and optimal matching problems'' funded by Università Italo Francese.  Research also partly funded by PNRR - M4C2 - Investimento 1.3, Partenariato Esteso PE00000013 - "FAIR - Future Artificial Intelligence Research" - Spoke 1 "Human-centered AI", funded by the European Commission under the NextGeneration EU programme.

%
%
%
%
%
%
%
%
%
%
%
%
%
%

\printbibliography

\end{document}